\magnification=\magstep1

\hsize=158truemm
\vsize=240truemm

\let\a=\ast
\let\p=\prime

\def\ltextindent#1{\hbox to \hangindent{#1\hss}\ignorespaces}
\def\lit{\par\noindent
               \hangindent=3\parindent\ltextindent}
\def\ref#1#2#3#4{\lit{[#1]} #2: {\it #3}, #4\medskip}
\def\reff#1#2#3{\lit{[#1]} #2: #3\medskip}
\def\llap#1{\hbox to 0pt{\hss$#1$}}
\def\rlap#1{\hbox to 0pt{$#1$\hss}}
\def\h#1{\hat #1}
\def\I{{I \kern-2pt I}}   
\def\cp{\hbox{$\times$\kern-1.8pt\vrule height
4.8pt depth 0pt width 0.3pt}\kern 1.2pt}

\def\bb#1{{\raise1.5pt\hbox{$^=$}}\kern-6.5pt #1}
\def\c#1{{\Cal{#1}}}
\def\h#1{\hat #1}
\def\o#1{\overline{#1}}
\def\rh#1{{\raise3.5pt\hbox{$\rightharpoonup$}}\kern-7.5pt #1}

\input amstex
\loadmsbm

\centerline{\bf On the topology of the group of invertible elements}

\bigskip
\centerline{ -- A Survey --}
\bigskip
\centerline{\it by Herbert Schr\"oder}
\bigskip
\bigskip\vskip 1cm

The topological structure of the group of invertible
elements in a unital Banach algebra (regular group for
short) has attracted topologists from the very beginning
of homotopy theory. Be it for its own sake simply to show
the instrumental power of newly invented methods or
because there were important applications, the most
notable one being perhaps the Atiyah-Singer index theorem
whose topological pillar is Bott's periodicity theorem
for the homotopy groups of the stable general linear
group. Recently, operator $K$-theory, which is the
homotopy theory of the stable regular group of a
$C^{\a}$-algebra, has been used to obtain index theorems in a
more general setting. While the properties that are
needed in index theory are by now quite well understood,
since there one only makes use of the stable regular
group, the topological structure of the regular group of a
$C^{\a}$-algebra itself is less well studied. Here we
want to survey the present state. Although today operator $K$-theory is merged
in $KK$-theory we do not include these new developments since they 
would take us to far and beyond our purpose.

The prototype of a regular group is the general linear
group $GL(n,\Bbb F)$ of invertible $n \times n$-matrices with
entries from $\Bbb F = \Bbb R$, $\Bbb C$, or $\Bbb H$. Exploring its
homotopical structure is intimately linked to the
development of algebraic topology, and even nowadays a
full understanding is out of reach. Without exposing the
highly sophisticated methods which are used to pursue
this problem, we will give a short historical account in the first section. 
In the second section we look at the general linear group
$GL(E)$ of invertible continuous linear operators on a
Banach space $E$, and also as an intermediate step from
finite to infinite dimensions at the Fredholm group
$GL_c(E)$, the subgroup of $GL(E)$ that consists of
perturbations of the identity by compact operators. 
General Banach algebras will be studied in section three, in particular,
commutative Banach algebras which for nearly forty years have taken an
independent development because of the close relation
to complex analysis. 
In the last section we deal exclusively with $C^{\a}$-algebras. We review the
important tools and results from operator $K$-theory and then go on
to discuss nonstable $K$-theory. However, we have to be selective
since operator $K$-theory is still rapidly expanding which makes it 
impossible to present all of the known $K$-groups. For more detailed
information the reader is referred to the
bibliography which is kept to a minimum and contains only 
articles that are well-suited for cross
reference so that most of the information can be retrieved from these.
In each section we also consider the topology of some associated spaces
such as the space of idempotents and the Grassmannian.
\plainfootnote{}{1991 AMS Subject Classification 46L\hfill Typeset
by \AmSTeX}

\vfill\break

\centerline{\bf 1. The classical groups} \bigskip

The simplest Banach algebras are the finite dimensional matrix algebras
$M(n,\Bbb F)$ (=$M(n,n,\Bbb F)$ if we denote by $M(n,m,\Bbb F)$ the space of $n \times m$-matrices
with real, complex or quaternionic entries). The regular group of $M(n,\Bbb F)$ is the
general linear group $GL(n,\Bbb F)$. The classical groups (a term coined by Hermann
Weyl) are subgroups of $GL(n,\Bbb F)$ and defined as follows. If $^-$ denotes
conjugation in $\Bbb F$ and $^t$ transposition in $M(n,m,\Bbb F)$, then
$$U(n,\Bbb F) = \{ X \in GL(n,\Bbb F) \mid \bar X^tX = I_n\}$$
($I_n$ the identity matrix), and the orthogonal groups, unitary groups, and
symplectic groups are given by $O(n) = U(n,\Bbb R)$, $U(n) = U(n,\Bbb C)$, and
$Sp(n) = U(n,\Bbb H)$, respectively. Accordingly, with
$$SU(n,\Bbb F) = \{ X \in U(n,\Bbb F) \mid \det X = 1\}$$
for $\Bbb F = \Bbb R$ and $\Bbb C$ one obtains the special orthogonal and the special unitary groups
$SO(n)$ and $SU(n)$, respectively, and with
$$Sp(n,\Bbb F) = \{ X \in GL(2n,\Bbb F) \mid X^tJ_nX = J_n\}$$
the real and the complex symplectic groups $Sp(n,\Bbb R)$ and $Sp(n,\Bbb C)$. Here $J_n$
denotes the matrix $\pmatrix0 & I_n\cr -I_n &0\cr\endpmatrix$. Moreover, one 
has the groups
$$U(n,m,\Bbb F) = \{ X \in GL(n+m,\Bbb F) \mid \bar X^tI_{n,m}X = I_{n,m}\}$$
where
$$I_{n,m} = \pmatrix
I_n&0\cr 0& -I_m\cr
\endpmatrix$$
and $\Bbb F = \Bbb R, \Bbb C,$ or $\Bbb H$.
Now polar decomposition provides a retraction from $GL(n,\Bbb F)$ onto $U(n,\Bbb F)$
which becomes a strong deformation retract via
$$(X,s) \mapsto X\cdot (\bar X^tX)^{-s/2}.$$
Restricting this homotopy to $Sp(n,\Bbb F)$ and to $U(n,m,\Bbb F)$ leads to the
following deformation retracts
$$Sp(n,\Bbb R) \cap O(2n) \cong U(n)$$
$$Sp(n,\Bbb C) \cap U(2n) \cong Sp(n)$$
$$U(n,m,\Bbb F) \cap U(n+m,\Bbb F) \cong U(n,\Bbb F) \times U(m,\Bbb F).$$
Finally, there are homeomorphisms $U(n) \to SU(n) \times S^1$ and
$O(n) \to SO(n) \times \Bbb Z_2$ which are both induced by the map
$$X = (X_1,\dots,X_n) \mapsto ((X_1\det X^{-1},X_2,\dots,X_n),\det X).$$
Therefore, it suffices to consider the compact groups $SO(n)$, $SU(n)$, and
$Sp(n)$.

The roots of algebraic topology are usually said to be found in
Euler's investigation of the K\"onigsberg bridges problem and
in his polyhedron formula. But he may also be considered to be the first
who has studied the classical groups when in 1770 he gave a parametrization 
of the rotation group $SO(n)$ by decomposing each transformation into a
product of two-dimensional rotations; of course, he did not yet conceive
$SO(n)$ as a topological group. In 1889 Kronecker showed that
$O(n)$ consists of two irreducible parts, the fact that constitutes chirality, 
i.e.\ the phenomenon of left-handedness and right-handedness in dimension $n=3$. He
also used the alternative parametrization
given by Cayley in 1846 to show that $SO(n)$ is an
$n(n-1)/2$-dimensional rational manifold. The rotation group
as a topological object appeared for the first time in 1897 in a
paper by Hurwitz [Htz] where he took up Euler's parametrization
to introduce volume elements on $SO(n)$ and on $SU(n)$. The point is a
footnote to this construction where he hints at the connectivity of these
groups. It was however up to Hermann Weyl to get fully aware of the real
importance that the type of connectivity has when one wants to deduce irreducible
representations of a Lie group from those of the corresponding Lie
algebra.

In 1924 Weyl announced that $SU(n)$ is simply connected and that
$SO(n)$ has a two-sheeted simply connected covering group if $n\ge 3$,
the so-called spin group $Spin(n)$. It is this double-covering 
$Spin(3)=SU(2)\to SO(3)$ that underlies the notion of spin in quantum theory
distinguishing between fermions and bosons. The complete proofs along with 
establishing
simple connectivity of $Sp(n)$ and finiteness of $\pi_1(G)$ for any Lie group
$G$ were presented in his famous series of
papers in 1925/6 [Wey1] -- that $\Bbb R$ covers $SO(2) = S^1$ via the
exponential map was already known when Poincar\'e introduced the
first homotopy or fundamental group $\pi_1(X)$ of a topological space
$X$ in a short note in 1892. We sketch Weyl's
argument for $SU(n)$: Given a closed path $C:[0,1] \to SU(n)$, $C(0) =
C(1) = I_n$, consider the ``spectral flow'' of the matrices $C(t)$, where
by a transversality argument if necessary one can assume that all
eigenvalues $\lambda_j(t)$, $j = 1,\dots,n$, $t\ne 0,1$, are simple. One
obtains $n$ closed paths $\lambda_j:[0,1] \to S^1$, which can be chosen
to satisfy $\sum \arg \lambda_i(t) = 2\pi$ and $\arg \lambda_1(t) < \cdots
<\arg \lambda_n(t)$ for $t \ne 0,1$. Thus these paths do not cross and 
are therefore homotopic to the constant path $\lambda(t) \equiv 1$.

Higher homotopy groups had been presented by E.\ \v Cech in
1931 at a meeting in Vienna and
again in 1932 at the ICM in Z\"urich but they did not receive much attention
until they were rediscovered in 1935 by Hurewicz [Hcz]. In the meantime
(in 1928) Elie Cartan [Ctn] had proved that the second homotopy group of a compact Lie
group $G$ vanishes. Still written in the language of homology theory
(vanishing of Betti numbers) his argument is similar to Weyl's noting
that any continuous map $S^2 \to G$ can be deformed to avoid the
singular set which has codimension 3 and, therefore, it is homotopic to a
constant one.

Further investigations of the homotopical structure of classical groups
have been intimately connected with the problem  of finding the maximal
number of linearly independent vector fields on spheres, in particular, to 
answer the question of parallelizability. In this context we have to mention
the work of Stiefel, Whitney, Eckmann, Ehresmann,
Feldbau, Hurewicz, and Steenrod on fibre bundles especially over spheres
and Stiefel manifolds continuing previous work of Freudenthal,
H.\ Hopf, and Pontrjagin on the topology of spheres.

The most important new device was the long exact homotopy sequence
which relates the homotopy groups of the total space $E$, of the base
space $B$, and of the fibre $F$ of a fibre bundle $(E,B,F)$:
$$\to \pi_k(F,\ast) \to \pi_k(E,\ast) \to \pi_k(B,\ast) \to
\pi_{k-1}(F,\ast) \to$$
(here $\a$ denote appropriate base points).
 When applied to the fibration
$$SU(n,\Bbb F) \to SU(n,\Bbb F)/SU(n-1,\Bbb F)\cong S^{n\cdot\dim\Bbb F-1}$$
it yields that the homotopy groups $\pi_k(SO(n))$, $\pi_k(SU(n))$, and
$\pi_k(Sp(n))$ do not depend on $n$ if $n \ge k+2$, $n\ge(k+1)/2$, and 
$n\ge(k-1)/4$, respectively, since $\pi_k(S^n) = 0$ if $0 \le k < n$ [Eck1].
We denote these stable homotopy groups by $\pi_k(SO)$, $\pi_k(SU)$, and 
$\pi_k(Sp)$.

Following Hurwicz' suggestion, Weyl used this
method in his book on classical groups [Wey2] to simplify his original
computation of the fundamental groups. At this occasion he also gives
the following nice argument that shows that the nontrivial loop
$$\phi(t) = \pmatrix
\cos t&-\sin t&0\cr
\sin t& \cos t&0\cr
0&0&1\cr
\endpmatrix,\quad 0\le t\le 2\pi,$$
in $SO(3)$ when cicumvented twice is homotopic to a constant one:
Consider two circular cones both of aperture $\alpha={\pi \over 2}$, one
being fixed in space the other rolling on the first. This movemant
describes a rotation by $4\pi$ about the axis of the moving cone. While 
increasing the angle $\alpha$
from ${\pi \over 2}$ to $\pi$ the motion of the rolling cone changes to a
mere wobble and, finally, at $\alpha = \pi$ comes completely to rest.

The homotopy groups $\pi_k(SO(n))$, $k=2,$ $3$, and $4$, have been computed
independently by Eckmann [Eck2] and G.W.\ Whitehead [Whi]
in 1941. Both also gave $\pi_5(SO(n))$ using a conjecture by Pontrjagin
who claimed that $\pi_5(S^3)$ is trivial. But it was only in 1950 that
this turned out to be wrong when Pontrjagin and Whitehead found that
$\pi_5(S^3) = \Bbb Z_2$. The corrected results were presented soon
afterwards in a survey article by Eckmann [Eck3] and in Steenrod's book
on the topology of fibre bundles [Ste]. Further computation of
$\pi_k(SO(n))$ has been done by Borel and Serre 
($k = 6$), Serre and Paechter ($k = 7$),
Sugawara ($k = 7$, 8, and 9, partial results for $k = 10$
and 11), and Toda [Tod].

However in several places the list given by Toda in 1955 did not agree
with results that Borel and Hirzebruch obtained in spring
1957 about the divisibility of certain characteristic classes [BH]. The
discussion around this controversy aroused Raoul Bott's interest, and in
1957 [Bot] he found the famous periodicity theorem and confirmed
Borel's and Hirzebruch's result.     \medskip

\lit{\bf Theorem } {\it The stable homotopy groups are periodic with
$$\pi_k(O) \cong \pi_{k+4}(Sp) \cong \pi_{k+8}(SO)$$ and
$$\pi_k(SU) \cong \pi_{k+2}(SU)$$ for $k\ge 0$.
}
\par\medskip

This theorem, a landmark in the history of homotopy theory, did not only
help to solve some of the classical problems mentioned before like the
parallelizability of spheres (according to Bott, Kervaire, and Milnor only
$S^1$, $S^3$ and $S^7$ have
trivial tangent bundle), it also provided the fundamentals of $K$-theory
which had been introduced by Atiyah and Hirzebruch. $K$-theory itself
found also many applications; we only mention the solution of the vector
field problem by Adams and the index theorem of Atiyah and Singer.

The periodicity theorem determined the homotopy groups of the classical
groups in the stable range. Afterwards many
topologists went on to explore the metastable range, i.e.\ $\pi_{2n+k}(SU(n))$,
$\pi_{n+k}(SO(n))$, and 
$\pi_{4n+k}(Sp(n))$ for $k$ small. A lot of computation has been done in the 
last thirty years above all by Japanese topologists. We do not want to go into 
details.
Instead, we refer to [Lun] which contains tables of the currently known
results as well as the credits.
We only mention that the methods are essentially to use
the exact homotopy sequence applied to fibrations of Stiefel manifolds. The
main problem here comes with the 2-primary components, i.e. the subgroups of
order a power of two. For a prime $p \ge 3$ the $p$-primary components
$\pi^p_k$ are much easier to find, for example for the rotation groups one 
has isomorphisms 
$$
\pi^p_k(Sp(n)) \cong \pi^p_k(SO(2n+1)),\quad n,k \ge 1.$$

For $n \ge 13$ and $k < n-2$ the homotopy groups $\pi_{n+k}(SO(n))$ can
be obtained from those of the Stiefel manifolds
$$
{St}_{n,m}(\Bbb F) = \{X \in M(n,m,\Bbb F) \mid \bar X^tX = I_m\}.$$
One usually denotes the real ones by $V_{n,m}$, the complex ones
by $W_{n,m}$, and the quaternionic ones by $X_{n,m}$. The crucial
relation found by Barrett and Mahowald in 1964 is then
$$
\pi_{n+k}(SO(n)) \cong \pi_{n+k}(SO) \oplus \pi_{n+k+1}(V_{2n,n}).$$

Homotopy groups $\pi_{k+p}(V_{k+m,m})$ are stable for
$m \ge p+2$. The metastable range has been explored by Nomura, see [Nom]
and references therein for earlier work. Complex Stiefel manifolds have
been inspected by Furukawa and Nomura [FN] and $X_{n,2}$ as the only
quaternionic Stiefel manifold considered so far by \^Oguchi [\^Ogu].

Finally, we want to mention the symmetric spaces $SO(2n)/U(n)$,
$SU(n)/SO(n)$, $Sp(n)/U(n)$, and $U(2n)/Sp(n)$, and the Grassmannians
$$
{Gr}_{n,m}(\Bbb F) = U(n,\Bbb F)/(U(m,\Bbb F) \times U(n-m,\Bbb F))$$
which have stable homotopy groups for $n$ large, too. These stable
groups also enter Bott's periodicity theorem. More precisely,
with the notation of the corresponding stable homotopy groups properly
understood one has
$$\eqalign{\pi_k(O) &\cong \pi_{k+1}(BO) \cong\pi_{k+2}(U/SO)\cong
\pi_{k+3}(Sp/U)\cr
\pi_{k}(Sp) &\cong \pi_{k+1}(BSp) \cong \pi_{k+2}(U/Sp)
\cong\pi_{k+3}(SO/U)\cr}$$
and $$\pi_k(SU) \cong \pi_{k+1}(BU),$$
where $BO$, $BU$, and $BSp$ are defined as $\lim {Gr}_{2n,n}(\Bbb F)$
with $\Bbb F = \Bbb R$, $\Bbb C$, and $\Bbb H$, respectively.
In the
metastable range the homotopy groups can be deduced from appropriate
homotopy sequences involving classical groups,
see again [Lun] for symmetric spaces; for Grassmannians
one uses the following relation due to James [Jam]
$$
\pi_k({Gr}_{n,m}(\Bbb F)) \cong \pi_k({St}_{n,m}(\Bbb F)) +
\pi_{k-1}(U(m,\Bbb F)),\quad 2m \le n.$$

\vfill\break

\centerline{\bf 2. The general linear group and the Fredholm group of a
Banach space} \bigskip

The possibility to distinguish between different orientations on
the Euclidean space $\Bbb R^n$ relies on the fact that the
orthogonal group $O(n)$ decomposes into exactly two components. To what extent
this carries over when one passes to an infinite dimensional real Hilbert space
$H_{\Bbb R}$ had been asked by Wintner in 1929 after he had shown that the
unitary group $U(H_{\Bbb C})$ of a complex Hilbert space $H_{\Bbb C}$ is connected -- by
the spectral theorem each unitary operator $U$ is of the form $U =
e^{iT}$ with
$T$ a self-adjoint operator hence $U(H_{\Bbb C})$ allows a
parametrization in the sense of Euler-Hurwitz [Win1]. Wintner gave an affirmative answer to
this question only about 20 years later in joint work with Putnam
[PW], although the prerequisites necessary for the proof were already
available in 1932 when Martin proved that any orthogonal
transformation $T$ can be written in the form $T = Q e^{-S}$ with $S$
skew-adjoint and $Q = I-2E$ a symmetry commuting with $S$.
Using this decomposition the proof is easy. Connect $e^{-S}$
with the identity $I$ by the path $t \mapsto
e^{-tS},\, 0 \le t \le 1$, and then, since $QE = - E$, $Q(I-E) = I-E$, and since
at least one of these projections is infinite, connect $E$ with $-E$ or
$I-E$ with $E-I$ according to as $E$ or $I-E$ is infinite. To do so and
finally to join $-I$ with $I$ one takes an orthonormal basis $(e_i)$ of the
corresponding infinite dimensional (sub)space (with obvious
modifications in the nonseparable case) and defines a family of orthogonal
transformations $U_t$ by
$$\eqalign{
U_t e_{2i-1} &= e_{2i-1}\, \cos t - e_{2i}\, \sin t\cr
U_t e_{2i}\, &= e_{2i-1}\, \sin t + e_{2i}\, \cos t\,.\cr}$$

In general $Q$ cannot be connected to $I$ by a path within the set of
self-adjoint orthogonal transformations (e.g. if $Q$ has positive and
negative eigenvalues) [Win2]. This is in sharp contrast to skew-adjoint
orthogonal transformations. Any two $S,T \in O(H_{\Bbb R})$ with $T^\a = - T$,
and $S^\a = - S$ belong to one and the same (connected) orbit, since the
orthogonal group acts transitively on the set of skew-adjoints by
conjugation [Win3].
Indeed, if one takes orthonormal bases $(e_i)$ and $(f_i)$
of $H_{\Bbb R}$ with
$$\eqalign{
&T\, e_{2i-1} = e_{2i} \,,\quad T\,e_{2i} =- e_{2i-1},\cr
&S\, f_{2i-1} = f_{2i} \,,\quad S\,f_{2i} = - f_{2i-1}\,,\cr}$$
then $Ue_i = f_i$ defines an orthogonal transformation with
$U^{-1} SU e_i = T e_i$.

Since any operator in $L(H_{\Bbb F})$, the Banach algebra of bounded linear
operators on the Hilbert space $H_{\Bbb F}$ ($\Bbb F = \Bbb R$ or $\Bbb C$), allows polar
decomposition, the same retraction and homotopy as in section 1 prove
that the general linear group $GL(H_{\Bbb F})$ of invertible
operators in $L(H_{\Bbb F})$ is connected.
This can also be seen in a
more elementary way using an idea of Nikolaas H. Kuiper's.
Given $T\in U(H_{\Bbb F})$ one successively constructs orthogonal unit vectors
 $a_1,a_2,\dots \in H_{\Bbb F}$ and two-dimensional
mutually orthogonal subspaces $A_1$, $A_2,\dots$ with $a_i,Ta_i \in A_i$. With $H^{\p}$
denoting the closure of the span of the $A_i$'s one defines by
 $$ S|H^{\p\bot} = I_{H^{\p\bot}}\,  \hbox{ and }\,
S|A_i\, \hbox{a rotation with}\, S(Ta_i) = a_i$$

\noindent an orthogonal transformation that can be connected with the
identity, i.e., $T$ can be connected with $ST$ and $ST$
is the identity on the subspace $H$ spanned by the $a_i$'s. 
With respect to $H_{\Bbb F} = H^{\bot}
\oplus H$ one gets 
$ST =
\pmatrix
T^{\p\p} & 0 \cr T^{\p} & I_H\cr
\endpmatrix$ and the path 
$\pmatrix
T^{\p\p} & 0 \cr tT^{\p} & I_H\cr
\endpmatrix$, $0 \le t \le 1$, leads to a
diagonal transformation. Now $H$ can be decomposed into countably many mutually
orthogonal subspaces isomorphic to $H^{\bot}$  and with respect to this
splitting one has
$$T^{\p\p} \oplus I_H =
T^{\p\p}\oplus T^{\p\p}T^{\p\p-1}\oplus I \oplus T^{\p\p}T^{\p\p-1}\oplus I
\cdots$$

\noindent Next one connects with
$$T^{\p\p}\oplus T^{\p\p-1}\oplus T^{\p\p}\oplus T^{\p\p-1}\oplus
T^{\p\p}\oplus\cdots$$

\noindent
and then with $I$ where in each step one uses the path
$$t \mapsto
\pmatrix
\cos t&-\sin t\cr \sin t &\cos t
\endpmatrix
\pmatrix
u&0\cr 0&I\cr
\endpmatrix
\pmatrix
\cos t&\sin t\cr -\sin t& \cos t\cr
\endpmatrix
\pmatrix
v &0\cr 0&I\cr
\endpmatrix,\, 0\le t \le {\pi\over 2},$$

\noindent
between
$$\pmatrix
uv&0\cr 0&I\cr
\endpmatrix \quad \hbox{ and }\quad 
\pmatrix
v&0\cr 0&u\cr
\endpmatrix.$$

In 1963 Shvarts and J\"anich noticed that
this trick can be used to contract the stable group 
$GL(\infty,\Bbb F) = \lim\limits_\to GL(n,\Bbb F)$ or even
the Fredholm group $GL_c(H_{\Bbb F}) =\hfill\break
\{T \in GL(H_{\Bbb F})\mid T-I\, \hbox{\rm compact}\}
$ to a point
within $GL(H_{\Bbb F})$.
Shvarts, Palais, and Atiyah then conjectured that already
$GL(H_{\Bbb F})$ should be contractible, and in due course a proof was
provided by Kuiper [Kui] using a refinement of
the previous idea. We cannot dwell on the consequences of this 
far-reaching result, e.g. it implies
that any Hilbert space bundle over a paracompact space is trivial and that the
unit sphere of Hilbert space is parallelizable.
Contractibility of $GL(H_{\Bbb F})$ in the strong topology had
been proved before by Douady and Dixmier [DD].

Of course, the question immediately came up to what extent Kuiper's result
carries over to arbitrary Banach spaces, since already Shvarts had
noticed that for some Banach spaces (among others $l^p,\, 1 \le p \le \infty$,
and $c_0$) the Fredholm group $GL_c(E)$ is contractible within
$GL(E)$ (see also [Deu2] for a recent variant). In 1965 Arlt proved that $c_0$,
the  space of sequences tending
to zero, is a ``Kuiper space'', but Douady found
examples of complex Banach spaces whose general linear groups are not even
connected. Specifically, he proved the following result, cf.\ the survey paper
[Mit]. 
\medskip
\lit{\bf Theorem 1} {\it If $E$ and $F$ are Banach spaces both isomorphic to
their hyperplanes and if $L(E,F) = L_c(E,F)$, then $GL(E \times F)$ has the homotopy
type of $$GL(E) \times GL(F) \times \Bbb Z \times  BGL(\Bbb F),$$
where $BGL(\Bbb F) = \lim\limits_{\to}GL(2n,\Bbb F)/(GL(n,\Bbb F) \times GL(n,\Bbb F))$.}
\par\medskip

Since then many Banach spaces have been inspected.
The classical sequence spaces $l^p$, $1 \le p < \infty$, (and again $c_0$)
were shown to be Kuiper spaces by Neubauer. On the other hand,
$l^p$ for $0 < p < 1$ is only a
metrizable topological vector space not even locally convex
but applying Neubauer's methods Cuellar [Cue] found that
the homotopy groups of the general
linear group also vanish while it is still not known whether it is contractible.
Neubauer's method has been refined by
Mitjagin and his coworkers leading to the following sufficient
criterion for contractibility.

\noindent
A Banach space $E$ is called weakly infinitely decomposable (WID) if
$E \cong E \times E$, if there
exists a total family of mutually orthogonal projections $P_k$ and
isomorphisms $T_k : P_kE \to E$, $k \ge 0$, such that
$$\eqalign{T_kP_kS & = T_{k+1}P_{k+1}, k \ge 0,\cr
T_kP_kS^{\p} & = T_{k-1}P_{k-1}, k \ge 1,\cr
P_0S^{\p} &= 0 \, \cr}$$
for some $S,S^{\p}\in L(E)$,
and if for any $B \in L(E)$ there exists ${\tilde B} \in L(E)$ with
$$ T_kP_k{\tilde B} = BT_kP_k, \quad k \ge 0.$$
$E$ has the property of smallness of operator blocks (SOB) if for any $\epsilon >
0$ and for any compact family of operators, $\{B\}$, there are orthogonal
projections $Q_1, Q_2$ with $Q_iE \cong E$, $i = 1,2$, and $\| Q_1 B Q_2\| \le
\epsilon$, for all $B$.

\medskip\lit{\bf Theorem 2}{\it If $E$ has properties WID and SOB then $GL(E)$
is contractible.}\par\medskip

\noindent
Using this general condition the following Banach spaces can be shown to be
 Kuiper spaces, cf. [Mit] for proofs or references:

\noindent
$C(K,F)$ where $F$ is a Banach space not having $c_0$ as a direct summand, and
where $K$ is either

\noindent -- an uncountable compact metric space,

\noindent -- an infinite compact topological group,

\noindent -- an infinite product of non-one-point compact metric spaces,

\noindent -- the Stone space of an infinite homogeneous measure algebra,

\noindent -- or the Stone-\v Cech compactification of the integers.

\noindent
This also holds if $K$ is the product of $\tau$ copies of a two-point space
with the topology induced by lexicographic order where $\tau \ge \omega$ is a
countable ordinal, cf. [Sem].

\noindent
Moreover, one has contractibility for 

\noindent -- $C^k(M)$, the space of k-times differentiable functions on a smooth
compact manifold $M$,

\noindent -- the Sobolev spaces $H^p_s(D)$, $1 < p < \infty,\, s \ge 0$, and
$W^p_\ell(D)$, $1 < p <\infty,\, \ell \in \Bbb N$, on a domain $D 
\subset \Bbb R^n$ with
regular boundary,

\noindent -- the spaces $L^p([0,1])$, $1 \le p \le \infty$, of $p$-integrable
measurable functions,

\noindent -- the spaces $L^p(\Omega, \mu)$, $1 < p < \infty$, of $p$-integrable
measurable functions on an arbitrary measure space $\Omega$,

\noindent
and for some classes of reflexive symmetric function spaces.

\noindent 
Another class of Banach spaces with contractible general linear group consists
of operator algebras, e.g. $L(H)$, $H$ a Hilbert space, or, more generally, 
$L(\c M)$, $\c M$ an injective von Neumann factor of infinite type or the 
hyperfinite factor of tye $II_1$ [SS], a countable decomposable von Neumann 
factor of type $III$ [Wil] and [SS], or $C_p(H)$, $1\le p \le \infty$, the von 
Neumann-Schatten classes.

There are also more non-Kuiper spaces. For example, Douady's theorem applies
to the following products $E\times F$ of Banach spaces:
$c_0 \times l^p$, $l^{p_1} \times l^{p_2},\, 1 \le p_1 < p_2 < \infty$, and
$J^{n_0}(p_0,\Bbb F) \times  J^{n_1}(p_1,\Bbb F)$, $1 \le p_0 < p_1 < \infty$,
where $J^m(p,\Bbb F)$ is the $m^{th}$ power of the
James space
$J(p,\Bbb F)$. More precisely,
$GL(J^m(p,\Bbb F))$ is homotopy
equivalent to $GL(m,\Bbb F)$ as is $GL(C(\Gamma_{m\omega_1},\Bbb F))$ where
$\Gamma_{m\omega_1}$ denotes
the ``long line'' of ordinals $\le m\omega_1$, $\omega_1$ the first uncountable
ordinal [Bel].
Finally, we mention that the general linear group $GL(X)$ of a real or
complex nuclear
Fr\'echet space $X$ with basis is dense in the algebra of bounded operators
$L_b(X)$ and is connected in the induced topology [Gra].

So far, we have only considered the general
linear group of a Banach space. Other ``classical'' groups can be
defined as well. In order to define orthogonal, unitary, or
symplectic groups one needs, respectively, a nondegenerate symmetric, 
hermitean, or skew-symmetric bilinear form $\varphi$ on a real (complex) 
Banach space $E$ which is then necessarily reflexive [Swa]. The subgroups of 
$GL(E)$ whose elements are $\varphi$-invariant are
real Banach Lie groups. Here we stick to Hilbert spaces, but see [Swa] for
symplectic groups on Banach spaces.
We have already met the orthogonal and unitary groups, $O(H_{\Bbb R})$ and
$U(H_{\Bbb C})$, which are deformation retracts of $GL(H_{\Bbb F})$, $\Bbb F = \Bbb R$ or
$\Bbb C$. The analogues of $U(n,m,\Bbb F)$ and $Sp(n,\Bbb R)$ are the group of
$Q$-unitary operators
$$
U(Q,H_{\Bbb F}) = \{T\in L(H_{\Bbb F})\mid T^\a QT = Q = TQT^\a\}\,,$$

\noindent
where $Q=Q^\a = Q^{-1} = 2P-1 \in GL(H_{\Bbb F})$, and the real symplectic
group
$$
Sp(J,H_{\Bbb R}) = \{T\in L(H_{\Bbb R}) \mid T^\a JT = J = TJT^\a\}\,,$$

\noindent
where $J\in O(H_{\Bbb R})$, $J^2 = -I$. 
The homotopy type of these groups is more or less folklore and is in
fact easily deduced from Kuiper's theorem. That  $U(Q,H_{\Bbb F})$
is connected has been noted in 1968 by Phillips [Phi], but was most
likely already known to Wintner in
the thirties. The same homotopy as in section 1 shows that $U(Q,H_{\Bbb F})$ 
retracts onto $U(PH_{\Bbb F}) \times U((1-P)H_{\Bbb F})$ [Ku\v c], and
that $Sp(J,H_{\Bbb R})$ is contractible [Swa].

When not
imposing an additional structure on a Banach space $E$ given by a bilinear form
one might think of the group of invertible
isometries $Iso(E)$ as the appropriate substitute for the group of unitaries.
But this group may have a rather peculiar topological structure and may even
fail to be a Banach Lie group [HK].
However, one has the result by Stern [Str] that there are real
Banach spaces whose group of invertible isometries is isomorphic to
any given group $G$ that contains a normal subgroup with two elements
$e$ and $\tau$. More
precisely, he proved that there exists a real Banach space $E$ and an
isomorphism $\Phi : E \to H$ with $\|\Phi\|\|\Phi^{-1}\| \le 1+\epsilon$
onto a Hilbert space $H$ such that $Iso(E) \cong G$ with $e$ corresponding
to $Id_H$ and $\tau$ to $-Id_H$.
In particular, $GL(E)$ is
contractible while $Iso(E)$ is any group extension by $\Bbb Z_2$. If $E$ is
a complex Banach space each component of $Iso(E)$ contains of course at least a
circle.

After this digression and before we come to the Fredholm group we take a brief
look at Stiefel and Grassmann manifolds in infinite dimensions. The
Grassmann manifold ${Gr}(E)$ of an infinite dimensional Banach space $E$ is
the set of complemented subspaces of $E$ with the opening topology (for not
necessarily complemented closed subspaces see [CrL]). The subset of complemented
subspaces isomorphic to a given complemented subspace $F \subset E$ is denoted
by ${Gr}_F(E)$. Specifically, we denote by ${Gr}_n(E)$ and ${Gr}_{\infty - n}(E)$ 
the set of $n$-dimensional and $n$-codimensional subspaces
of $E$, respectively. The set ${St}_F(E)$ of
monomorphisms of a complemented subspace $F$ into $E$ with complemented range
generalizes the Stiefel manifold in infinite dimensions. $GL(F)$ acts on  
${St}_F(E)$ and ${Gr}(E) = {St}_F(E)/GL(F)$. If $F$ is an $n$-dimensional
subspace of $E$ then ${St}_F(E)$ is contractible [DD], and hence the exact
homotopy sequence of fibre bundles determines the homotopy type of 
${Gr}_n(E)$ and of ${Gr}_{\infty - n}(E)$. In case of a Hilbert space $E$
the various homotopy types are completely known, cf. [Luf1]; in
general, they depend on the homotopy type of $GL(E)$.

Whereas $GL_c(E)$ is contractible within $GL(E)$ for some Banach spaces $E$
the Fredholm group $GL_c(E)$ itself is homotopy equivalent to $GL(\infty,\Bbb F)$.
This was shown by Shvarts and Palais in 1963 for a Hilbert space and
in general independently by Elworthy [ET] and
Geba [Geb] in 1968. Already Palais had considered the groups
$GL_{C_p}(H_{\Bbb F})$ defined analogously using the symmetric ideals $C_p(H_{\Bbb F})$
 of von Neumann-Schatten operators instead of the ideal of compacts $C(H_{\Bbb F})$.
These groups are again homotopy equivalent to $GL(\infty,\Bbb F)$.
Geba's proof implies the following more general result stated in [dlH1]:
\medskip\lit{\bf Theorem 3} {\it If $E$ is a Banach space and $P(E) \subset L(E)$
a normed (not necessarily closed) subspace with}
\lit{} {\it i) $1+X\in {Fred}(E)$ (= set of Fredholm operators)}
\lit{} {\it ii) $X+C_0(E) \subset P(E)$ ($C_0(E)$ the ideal of finite rank
operators)
 for all $X\in P(E)$ and}
\lit{}{\it iii) $C_0(E) \times P(E) \ni
(T,X) \mapsto TX \in P(E)$ continuous,}
\lit{} {\it then $GL_p(E) = \{1+X\in GL(E)\mid X\in P(E)\}$ is homotopy
equivalent to $GL(\infty,\Bbb F)$.}\par\medskip 

In particular, $GL_{C_p}(H_{\Bbb R})$ and its deformation retract
$O_{C_p}(H_{\Bbb R})$ have two connected components and the component
$SO_{C_p}(H_{\Bbb R})\subset O_{C_p}(H_{\Bbb R})$ that contains the identity has a two-sheeted universal
covering group $Spin_{C_p}(H_{\Bbb R})$, which can be realized as a
subgroup of the unitary group of the infinite dimensional Clifford
$C^\a$-algebra in case that $p = 1$ [dlH2], and of the hyperfinite
$W^\a$-algebra factor of type $II_1$ in case that $p=2$ [Ply]. 

Moreover, for $J\in O(H_{\Bbb R})$ with $J^2 = -I$, $J-J^\a\in \c S$ one
defines  $$
G_{C_p}(H_{\Bbb R},J) = \{T\in GL(H_{\Bbb R}) \mid TJ - JT \in C_p(H_{\Bbb R})\}$$

\noindent
with retract
$$
O_{C_p}(H_{\Bbb R},J) = G_{C_p}(H_{\Bbb R},J) \cap O(H_{\Bbb R}) \cong
O/U,$$

\noindent
and for $Q\in U(H_{\Bbb F})$, $Q=Q^\a$, the group
$$
G_{C_p}(H_{\Bbb F},Q) = \{T \in GL(H_{\Bbb F}) \mid TQ - QT\in C_p(H_{\Bbb F})\}.$$

\noindent
They appear as regular groups of the Banach algebras of operators $T \in L(H_{\Bbb F})$
with $$\|T\|_{p,V} = \|T\| + \|TV - VT\|_{\c S} < \infty,$$
where $V = J$ or $Q$, cf.\ [CE] and previous work by Carey and coauthors cited there. 
$G_{C_1}(H_{\Bbb C},Q)$ was
denoted by $G_{res}(H)$ in [PrS], cf.\ also [Woj].
More generally, Carey and Evans [CE] considered Banach algebras
consisting of operators that commute with up to one generator of a Clifford
algebra $\c C_k$, the commutator
with the last generator lying in a symmetric ideal $\c S$ of compact operators.
The periodicity of the Clifford algebras (of period 8 in the real and of period
2 in the complex case) is inherited by these algebras
and the principal components of the corresponding regular groups have
the same homotopy type as the stablilized homogeneous spaces encountered
at the end of section 1. Moreover, suitably
modified they can be used as classifying spaces in $KK$-theory.

The classical example of a classifying space is the set ${Fred} (H_{\Bbb F})$ 
of Fredholm operators which is the preimage of the regular group of the
Calkin algebra $Cal(H_{\Bbb F}) = L(H_{\Bbb F}) / C(H_{\Bbb F})$  
under the canonical quotient map $\pi$. Here $\pi$ induces even a homotopy
equivalence between ${Fred} (H_{\Bbb F})$ and $G(Cal(H_{\Bbb F}))$, cf. the
appendix in [Ati]. The components of ${Fred} (H_{\Bbb F})$ are separated by
the index map which indeed is an isomorphism between the set of components
and $\Bbb Z$.
This isomorphism has been generalized by Atiyah and Singer
[AS] and by Karoubi [Kar2] as follows:
\medskip
\lit{\bf
Theorem 4}{\it For any compact Hausdorff space $X$, and $\Bbb F = \Bbb R, \Bbb C$, or
$\Bbb H$ one has $$K^n_{\Bbb F}(X) \cong [X,F^n(H_{\Bbb F})],$$ where
$F^0(H_{\Bbb F}) = {Fred} (H_{\Bbb F})$ and where  $F^n(H_{\Bbb F})$ is a
subset of ${Fred} (H_{\Bbb F})$ characterized by certain spectral
properties.}
\par\medskip

\noindent
For equivariant $K$-theory $K_G^n(X)$, $G$ a compact Hausdorff group, one
has an isomorphism of $K_G^n(X)$ with $[X,F^n(H_{\Bbb F})]_G$ when $H_{\Bbb F}$ is a
$G$-module and $F^n(H_{\Bbb F})$ is considered as a $G$-space ($[\cdot,\cdot ]_G$
denotes $G$-homotopy classes of $G$-invariant maps). This is proved for
$n = 0$ by Segal and
Matumoto and in general by Ku\v cment and Pankov [KP] using an
equivariant version of Kuiper's theorem given by Segal. Note that in these
papers it is tacitly assumed that the $G$ action is norm continuous.
While the equivariant contractibility is by the time still open the above
isomorphism nonetheless holds, cf.\ [Phs2].

There are partial extensions to Fredholm operators (or semi-Fredholm
operators) on a Banach space $E$, see [ZKKP] for the following and for
further references. Here one has an exact sequence
$$[X,GL_c(E)] \to [X,GL(E)] \to [X,{Fred}(E)] \to K_{\Bbb F}^0(X),$$
valid for any compact Hausdorff space $X$. If $E$ contains a complemented
infinite dimensional subspace with symmetric basis then this leads to the exact
sequence
$$0 \to [X,GL(E)] \to [X,{Fred}(E)] \to K^0_{\Bbb F}(X)\to 0$$
and in case that $E$ is a Kuiper space to the isomorphism 
$$K^0_{\Bbb F}(X) \cong [X,{Fred}(E)].$$

If $E$ and $F$ are infinite dimensional Banach spaces and if the set 
${Fred}(E,F)$ of Fredholm operators in $L(E,F)$ is nonempty then $GL(E)$ is
connected (contractible) if and only if $GL(F)$ is. Moreover, in
this case the components of ${Fred}(E,F)$ consist just of Fredholm operators
with the same index [GK].

Also the Calkin algebra has been inspected in more detail. The homotopy type of
the regular group is completely determined by the previous theorem. The set of
projections consists of three arcwise connected components, $\{0\}$, $\{1\}$,
and the set of nontrivial projections. Each component is simply connected, and
so is each component of the regular group [KK].

\vfill\break

\centerline{\bf 3. General Banach algebras}
\bigskip

Let $\c B$ be a real or complex Banach algebra with unit and let $G(\c B)$
denote the group of invertible elements. $G(\c B)$ is a Banach
Lie group and an open subset of $\c B$. Therefore, $G(\c B)$ is locally arcwise
connected and its components are arcwise connected. We denote the principal
component, that is the one which contains the unit element, by $G^0(\c B)$.
$G^0(\c B)$ is a Banach Lie group, too, and a normal subgroup of $G(\c B)$; the
quotient group $I(\c B) = G(\c B)/G^0(\c B)$ is called the index group of $\c
B$. Besides the regular group $G(\c B)$ we also consider the set 
$${Id}(\c B) = \{ p \in \c B  \mid  p^2 = p\}$$ 
of all idempotent elements
which is homeomorphic to the set of symmetries
$${Sym}(\c B) = \{ x \in \c B  \mid  x^2 = 1\}.$$
If $\c B$ is real one also has 
$${D}(\c B) = \{ x \in \c B  \mid  x^2 = -1\}.$$

\noindent
It has been shown by various authors that
${Id}(\c B)$ is locally arcwise connected in the
relative topology and that the components are arcwise connected. More precisely,
 the components are equal to the orbits of $G^0(\c B)$
acting on ${Id}(\c B)$ by conjugation, cf.\ the references in [CPR]. 
The same holds for the action
of $G^0(\c B)$ on ${D}(\c B)$ [Fuj2].
Corach, Porta, and Recht [CPR] give a more general result: If $P(x)$ is a
polynomial with simple roots then
$$A = \{ a \in \c B \mid  P(a) = 0\}$$ 
is locally arcwise connected
and the component containing $a$ is just again the orbit
$\{gag^{-1} \mid g \in G^0(\c B)\}$.

On ${Id}(\c B)$ one defines an equivalence relation by 
$$ p \sim  q  \quad \hbox{iff} \quad pq	 = q \quad \hbox{and}\quad qp = p.$$
The quotient space ${Gr}(\c B) = {Id}(\c B)/_{\sim} $ is called the
Grassmannian of $\c B$. The main result here is [PR1]:

\medskip\lit{\bf Theorem 1}{\it ${Gr}(\c B)$ and ${Id}(\c B)$ are
homotopy equivalent.}\par\medskip

\noindent
Because given $p \in {Id}(\c B)$ one can find a neighborhood $U$ with $1 -
(q-p)^2$ invertible for $q \in U$ and $q \mapsto [1 - (q-p)^2]^{-1}q$
continuous on $U$. Then $h : U \ni q \mapsto [1 - (q-p)^2]^{-1}qp \in 
{Id}(\c B)$ is a continuous map. Using a partition of unity one defines $h$
globally. Now $h : {Id}(\c B) \to {Id}(\c B)$ is continuous and has the
following properties
$$\eqalign{ph(p) &\sim p \quad \hbox{for all}\quad p, \cr
h(p) &= h(q) \quad\hbox{iff}\quad p \sim q.\cr}$$
The same proof carries over to topological algebras with open regular group,
continuous inversion, and ${Gr}(\c B)$ paracompact [PR2], cf. also the
paper [Grm] by Gramsch where these results are obtained by explicit
calculations using ``rational coordinates''.
If $\c B$ is a
$C^{\a}$-algebra, ${Gr}(\c B)$ can be identified with the set of
projections, i.e. self-adjoint idempotents.

Abstract Banach algebras had been introduced in 1936 by Nagumo as
 linear metric rings. In 1939 Gelfand founded the theory of
maximal ideals for commutative Banach algebras and showed that the set of all
maximal ideals can be endowed with a topology making it into a locally compact
Hausdorff space which is compact if the algebra has a unit. If we identify
maximal ideals with kernels of characters, i.e., nontrivial homomorphisms
$\varphi : \c B \to \Bbb C$, this topology is nothing but the weak-$*$-topology on the
set of characters. Via the correspondence between the space of maximal ideals
$M(\c B)$ and the set of characters $\Omega(\c B)$ each $x \in \c B$ defines a
continuous function $\hat x \in C(M(\c B))$ by $\hat x(\omega) = \omega(x)$,
$\omega \in \Omega(\c B)$. The assignment
$$\c B \ni x \mapsto \hat x \in C(M(\c B)),$$
the so-called Gelfand transform, is a homomorphism which in general is neither
injective nor surjective. 

Even in the commutative case investigations of the homotopy type of the regular
group of a Banach algebra have not been taken further than to determine the
index group and the fundamental group of the principal component.
The first algebra under inspection was the algebra $CAP(\Bbb R)$ of continuous 
almost periodic functions on $\Bbb R$ (here $M(CAR(\Bbb R)) = \Bbb R^B$, the so-called Bohr
compactification of $\Bbb R$). In 1930 H.\ Bohr proved a conjecture by
Wintner who had claimed that the regular group of $CAP(\Bbb R)$ has
uncountably many components labeled by the ``mean motion'' or ``mean
winding number''.  For $f \in CAP(\Bbb R)$ this winding number is defined by
$$\tau(f) = \lim_{T \to \infty} {1 \over 2T} \int_{-T}^T d({arg} \, f)$$
and generalizes the usual winding number $\tau(f)$ of a periodic
function $f$ defined by
$$\tau(f)= {1 \over 2\pi i} \int_{S^1} f^{-1} df = {1
\over 2\pi} \int_0^{2\pi} d({arg}\, f).$$ 
Both maps induce isomorphisms $I(G(CAP(\Bbb R))) \cong \Bbb R$ and $I(C(S^1)) \cong
\Bbb Z$.
In 1936 van Kampen has extended this result to almost periodic
functions on a connected Lie group $G$,
proving that to each connected component of the regular group $G(CAP(G))$ there
corresponds a group character.
Since $CAP(S^1) = C(S^1)$, one recovers once again the classical result
$I(C(S^1)) \cong \Bbb Z$.

For a compact metric space $X$ Bruschlinsky had already shown in 1934
that the index group of $C(X)$ coincides with its first \v Cech
cohomology group $H^1(X,\Bbb Z)$, and so contains in particular no element
of finite torsion, a fact established by Lorch in 1942 for any commuative
Banach algebra.
That $I(\c B)$ is commutative for any commutative Banach algebra has 
been proved independently in 1960 by
Royden and Arens; more precisely:

\medskip 
\lit{\bf Theorem 2}{\it For a commutative Banach algebra $\c B$ with unit one 
has  
$$I(\c B) \cong H^1(M(\c B),\Bbb Z).$$ }\par\medskip

\noindent
For noncommutative Banach algebras this turns out to be wrong in general, $I(\c
B)$ may be neither commutative nor torsion free. To see this consider the index
group of $\c B = C(S^k,M(n,\Bbb C))$ which is isomorphic to $[S^k,U(n)]$, and hence
$$I(\c B) \cong \pi_k(U(n))$$
since $U(n)$ is arcwise connected. This gives finite nontrivial index groups
for suitable $n$ and $k$, cf.\ the references in section 1. A simpler example of a finite (commutative) index group has been
found by Paulsen [Pau]:
$$SM(n,\Bbb C) = \{ f \in C([0,1],M(n,\Bbb C)) \mid f(0), f(1) \in \Bbb C I_n\},$$
the so-called unreduced suspension of $M(n,\Bbb C)$, has index group isomorphic
to $\Bbb Z_n$.
A nonabelian index group is obtained for $\c B =
C(U(2) \times U(2),M(2,\Bbb C))$. Here
$$I(\c B) \cong [U(2) \times U(2),U(2)]$$
is nonabelian since according to I.M.\ James the multiplication in $U(2)$ 
is not homotopy commutative, i.e.
$$(X,Y) \mapsto XY  \quad \hbox{and} \quad (X,Y) \mapsto YX$$
are not homotopic [Yue].
 
It is an open problem to determine the groups that can be realized as the index
group of some Banach algebra or that arise as index group of the operator
algebra $L(E)$ for some Banach space $E$.

The fundamental group of the principal component of a commutative Banach
algebra $\c B$ is isomorphic to the kernel of the exponential map $\exp :
\c B \to G^0(\c B)$ which is onto and induces a covering map as in the classical
case $\c B = \Bbb C$. Combining this observation by Blum from
1952 [Blu] with Shilov's decomposition theorem that
characterizes the idempotents of $\c B$ by the closed and open sets
in $M(\c B)$ one obtains

\medskip
\lit{\bf Theorem 3} {\it The fundamental group $\pi_1(G^0(\c B))$ of the
principal component $G^0(\c B)$ of a unital commutative  Banach algebra $\c B$ is
isomorphic to $H^0(M(\c B),\Bbb Z)$.}\par\medskip

So far we have assumed without saying that $\c B$ is a complex commutative
Banach algebra. Indeed, the structure of real commutative Banach algebras is
more complicated. For example, a real commutative $C^{\a}$-algebra is in
general not of the form $C(X,\Bbb R)$,
but isomorphic to  
$$ C(X,\tau) = \{ f \in C(X) \mid \o{f(x)} = f(\tau x)\,,\; x\in X\}\,,$$

\noindent
where $\tau$ is an involutive homeomorphismus of $X$. Accordingly, the real
analogues of the Arens-Royden theorem proved by Alling and
Campbell [AC] in 1971 and by Furutani [Frt] in 1975 are
more involved. 
\medskip
\lit{\bf Theorem 4} {\it Let $\c B^\sharp$ denote the set of all $a \in G(\c
B)$ such that to each $N \in M(\c B)$ with $\c B/N \cong \Bbb R$ there exists $a_N >
0$ with $a - a_N\in N$. Then $\c B^\sharp / \exp \c B$ is isomorphic to
$H^1(M(\c B),\Bbb Z^G)$, the first cohomology group of $M(\c B)$ with coefficients
in a certain sheaf of abelian groups $\Bbb Z^G$.}\par\medskip

\noindent
The simpler results of Lorch and Blum had already been carried over by
Ingelstam [Ing] in 1964. For a real commutative Banach algebra $\c B$,
not containing a complex subalgebra, he proved that $I(\c B)$ is either infinite or
isomorphic to the group of idempotents and that $G^0(\c B)$ is simply
connected in this case.

For quite a time it has been tried to bring also the higher \v Cech homology
groups into the game, but successfully so far only for $k \le 3$ [Tay3].
A different direction has been taken by Arens in 1965 when he
looked at the regular group $GL(n,\c B)$ of the tensor product
$\c B\otimes M(n,\Bbb C) = M(n,\c B)$, cf.\ the survey
articles by Forster [For] and by Taylor [Tay1].
Using a deep theorem of Grauert from complex analysis he
found that 
$$
[GL(n,\c B)] \cong [M(\c B),GL(n,\Bbb C)]\,.$$

\noindent
The relation between $GL(n,\c B)$ and $C(M(\c B),GL(n,\Bbb C))$
has been studied further by Eidlin and Novodvorskii.
In particular, Novodvorsikii's methods have been exploited later by Taylor
[Tay2] and Raeburn [Rae] who proved the following result.
\medskip
\lit{\bf Theorem 5} {\it Let $F$ be a closed Banach submanifold of an open subset
in a Banach space $E$. If $F$ is a discrete union of Banach homogeneous spaces
then
$$ \c B_F := \{a\in \c B\h\otimes E \mid a = Sf\,,\; f\in \c O(U,F)\,,\;
U\supset M(\c B) \}$$ (where $S$ is the extension of functional calculus
and $\c O(U,F)$ is the set of analytic functions from $U$ to $F$) 
is locally arcwise connected and the 
Gelfand transform induces a bijection 
$$
[\c B_F] \to [M(\c B),F]\,.$$}\par\medskip

\noindent
All of the previous results can be obtained by specialization, e.g.\ Arens'
result by taking  $F=GL(n,\Bbb C) \subset \Bbb C^{n^2}$. For
$E=\Bbb C^n$ the theorem had been proved in [Tay2]. In this special case $\c B_F$
can be replaced by the  homotopy equivalent set
$$
\c B^F = \{ a = (a_1,\cdots,a_n) \in \c B^n \mid \sigma (a) \subset F\}$$

\noindent
where $\sigma(a)$ denotes the Taylor spectrum of $a$ [LZ]. Several
corollaries are stated in [Tay2], [Rae], and in papers by
Fujii [Fuj1-3].
For example, if $E = \c A$ is a unital Banach algebra and $F=G(\c A)$ Raeburn
recovers a theorem by Davie which even
holds for Fr\'echet algebras.
\medskip
\lit{\bf Corollary 1} {\it If $\c A$ is a unital Banach algebra then the Gelfand
transform induces a homotopy equivalence between $G(\c B\h\otimes \c A)$ and
$C(M(\c B),G(\c A))$.}\par\medskip

\noindent
Analogous results can be obtained taking $F= {Id}(\c A)$ [Rae], 
${D}(\c A)$  [Fuj2], or ${Gr}(\c A)$ [PR1].

Since we are mostly interested in higher homotopy groups $\pi_k(G^0(\c B))$
we take $E = C(S^k)$ and get a homotopy equivalence between
$$
G(\c B\h\otimes C(S^k)) = G(C(S^k,\c B)) = C(S^k,G(\c B))$$
and
$$\eqalign{C(M(\c B),G(C(S^k)))
&= G(C(M(\c B),C(S^k)))\cr
&= G(C(M(\c B)\times S^k))\cr
&= C(M(\c B) \times S^k,\Bbb C \setminus \{0\})\,,\cr}$$

\noindent
Thus the computation of $\pi_k(G^0(A))$ is reduced to the computation of the
cohomotopy groups $\pi^1(M(\c B)\times S^k)$  -- a by no means easier
problem. It seems to become even more difficult if one takes $E=C(S^k,M(n,\Bbb C))$
instead.
However, passing to the direct limit $\lim\limits_{n\to \infty} C(S^k,M(n,\Bbb C))$
leads into the field of $K$-theory and simplifies in the following result
first obtained by Novodvorskii and later by Taylor.

\medskip
\lit{\bf Corollary 2} {\it For any complex commutative Banach algebra $\c B$ one has}
$$K_i(\c B) \cong K^{-i}(M(\c B)),\, i = 0,1.$$
\medskip

\noindent
Here $K_i(\c B)$ is the $K$-theory of a Banach algebra defined as classes of 
stable equivalent finitely generated projective $\c B$-modules if $i = 0$
and as $\pi_0(GL(\infty,\c B))$ if $i = 1$, cf.\ [Tay1], and $K^{-i}(M(\c B))$ is
the topological $K$-theory of the compact space $M(\c B)$, cf.\ [Ati] and [Kar3]. In
particular, the Chern character of $K$-theory and ``rationalization''
lead to a weak characterization envisaged earlier involving higher \v Cech
cohomology groups, since
$$
K^{-i} (M(\c B)) \otimes \Bbb Q \cong \bigoplus_{k\ge 0} H^{2k+i} (M(\c
B),\Bbb Q)\,.$$
The $K$-groups of some topological spaces are known but we can neither
discuss the techniques to obtain them nor display the results since this
would be beyond the scope of our survey (some special examples will
appear in the next section).

$K$-theory for Banach algebras is not restricted to the commuative case.
Most of the theory in [Tay1] is also valid for general Banach algebras (even
graded [vDa] and [Kan]) or Fr\'echet algebras [Phs1].
The fundamental result is the generalized periodicity theorem proved by
Wood in 1965 [Woo], cf. also [Kar1,3].
\medskip
\lit{\bf Theorem 6}{\it If $\c B$ is a unital Banach algebra then
$$K_0(\c B) \cong \pi_{2k-1}(GL^0(\infty,\c B))\, \hbox{and}\,
K_1(\c B) \cong \pi_{2k}(GL^0(\infty,\c B)),\, k\ge 1,$$
where $GL^0(\infty,\c B)$ denotes the principal component of $GL(\infty,\c
B)$.}\par\medskip
There are more parallels to the homotopical structure of the classical
groups, at least in special cases: if the Bass stable rank $sr(\c B)$ of
$\c B$ is finite according to Corach and Larotonda (see [CoL],
[Cor], and [CoS2]) the homotopy groups $\pi_k(GL(n,\c B))$ stabilize for
large $n$. Before we state the proper result, recall that the Bass stable rank
$sr(\c B)$ is the smallest number $n$ (or $\infty$ if no such $n$
exists) such that for all $(a_1,\cdots,a_{n+1}) \in \c B^{n+1}$ with 
$\sum a_i \c B= \c B$ there is a $(b_1,\cdots,b_n) \in \c B^n$ with 
$
\sum (a_i + b_i a_{n+1}) \c B = \c B$.
For  $\c B = C(X,\Bbb F)$ one has $sr(\c B) = [{1\over{\dim \Bbb F}}\, \dim  X] + 1$
and if $\c B$ is a $C^\a$-algebra then $sr(\c B)$ coincides with the
topological stable rank $tsr(\c B)$ introduced by Rieffel [Rie1].
$tsr(\c B)$ is the smallest number $n$ such that the set of $n$-tuples of
elements of $\c B$ which generate $\c B$ as a left (right) ideal is dense
in $\c B^n$ or in case of a $C^\a$-algebra that $R_n(\c B) =  \{x\in \c
B^n\mid \sum x_i^\a x_i \in G(\c B)\}$ is dense in $\c B^n$.
Note that $sr(\c B)$ and $tsr(\c B)$ may differ for a commutative Banach
algebra $\c B$, e.g. the polydisc-algebra $A(\Bbb D^n)$ has stable rank
$[{n \over 2}] + 1$ and topological stable rank $n+1$, cf. [CoS1].
\medskip
\lit{\bf Theorem 7} {\it If $\c B$ is a unital
Banach algebra (or a topological algebra with open
regular group) with finite stable rank then
$$\pi_k(GL(n,\c B)) \cong \pi_k(GL(n+1,\c B))$$
and
$$\pi_k({Id}(M(n,\c B))) \cong \pi_k({Id}(M(n+1,\c B)))$$  
if  $n\ge sr(\c B) + k+1$ .}\par\medskip

\noindent
If $sr(\c B)$ is not finite both assertions may be wrong for any $n$.
Already the index group  $I(M(n,\c B)) = \pi_0(GL(n,\c B))$ may be unstable in
this case, cf.\  [Srd4].

Finally, we want to look at Banach
algebras $\c A$ and $\c B$ which are related by a
dense injective morphism $\varphi :\c A \to \c B$ and see how this influences
the homotopy type of the corresponding regular groups. First we have
Karoubi's density theorem [Kar4] (which in a special case
goes back to Atiyah [Ati]).

\medskip
\lit{\bf Theorem 8}{\it If $\varphi : \c A \to \c B$ is a dense injective morphism
of unital Banach algebras, i.e., $\varphi(\c A)$ is dense in $\c B$, and if
$$GL(n,\c B) \cap M(n,\c A) = GL(n,\c A)$$
for any $n \ge 1$ (identifying $\c A$ with its image), then $\varphi$ induces
isomorphisms
$$\varphi_{\a} : K_i(\c A) \to K_i(\c B),\, i \ge 0.$$}\par\medskip

\noindent
In particular, this is guaranteed for complex Banach algebras if $\c A$ is stable
under holomorphic functional calculus in $\c B$, i.e., if for any $n$ and $a \in
M(n,\c A) \subset M(n,\c B)$ one has $f(a) \in M(n,\c A)$ for $f$ holomorphic in
a neighborhood of the spectrum of $a$ in $M(n,\c B)$, cf.\ [Bos].

A stronger result has been proved by Bost [Bos]. As a
consequence of the following theorem one finds that the conclusion of
Theorem 8 still holds if only $G(\c B) \cap \c A = G(\c A)$ is satisfied.
\medskip
\lit{\bf Theorem 9}{\it If $\c A$ and $\c B$ are unital Banach algebras (or
Fr\'echet algebras with open regular group and continuous inversion) and if
$\varphi : \c A \to \c B$ is a continuous morphism with $\varphi(\c A)$ dense
in $\c B$ and $\varphi^{-1}(G(\c B)) = G(\c A)$ then the induced maps
$$\varphi_n : M(n,\c A) \ni (a_{ij}) \mapsto (\varphi(a_{ij})) \in M(n,\c B)$$
define homotopy equivalences
$$\varphi_n : {Id}(M(n,\c A))  \to {Id}(M(n,\c B))$$ 
and 
$$\varphi_n : GL(n,\c A) \to GL(n,\c B)$$
for $n \ge 1$.}\par\medskip

\noindent
Bost has also considered many special cases including subalgebras of smooth
elements with respect to a group action or crossed products by abelian
groups.

\vfill\break

\centerline{\bf 4. The regular group of a $C^\a$-algebra} 
\bigskip

In his book on Banach algebras [Ric] Rickart shows that
Wintner's result about the connectivity of the unitary group carries over when
the operator algebra $L(H_{\Bbb C})$ is replaced by a Rickart $C^\a$-algebra. Indeed,
the only condition to meet is that all maximal commutative $C^\a$-subalgebras
 are of the form $C(X)$ with totally disconnected $X$.
In particular, the unitary group of an $AW^\a$-algebra and hence of a
$W^{\a}$-algebra are connected. If $\c M$ is a real $W^{\a}$-algebra without
finite discrete part essentially the
same arguments that were used by Putnam and Wintner show that
its orthogonal group is also connected. This has been discovered by
Ekman [Ekm] and Schreiber (unpublished) in 1974.
Even Kuiper's theorem appears as a special case of a more general result
about $W^{\a}$-algebras. After a series of papers by Breuer, Singer, and
Willgerodt dealing with countably decomposable ones,
Br\"uning and Willgerodt [BW] finally proved contractibility
for the regular group of any real or complex
$W^\a$-algebra of type $I_\infty$, $II_\infty$, or\ $III$.

\medskip\lit{\bf Theorem 1}{\it If $\c M$ is a $W^{\a}$-algebra of infinite type
then $G(\c M)$ and hence $U(\c M)$ are contractible.}\par\medskip

A finite continuous $W^\a$-algebra behaves
differently. Araki, L.\ Smith, and M.-S. Bae Smith [ASS] observed
that the fundamental group $\pi_1(G\c M)$ of a type $II_1$-factor $\c M$ is
nontrivial, and Handelman proved that this extends to the nonfactorial case [Han1].
More precisely, they proved that
$$\pi_1(G(\c M)) \cong K_0(\c M) \cong C(\Omega,\Bbb C),$$
if $\Omega$ is the Stonean space of the center of $\c M$.
Handelman also considered the fundamental group of the regular group of
an $AF\, C^{\a}$-al\-ge\-bra; the special case of $UHF$-algebras had been treated
before by Singer, L.\ and M.-S.\ Bae Smith in unpublished work. In the
factorial case (and also in case of a $UHF$-algebra) $K_0(\c M)$ is isomorphic to
$\Bbb R$ (or to a dense subgroup of $\Bbb R$) and the isomorphism is induced by
assigning the winding number $\tau(f) = {1 \over 2\pi i} \int_{S^1} tr
(f^{-1}df)$ to a differentiable loop $f$. Since the dimension group is generated
by the traces of projections, taking simple loops of the form $t \mapsto e^{2\pi
i pt}$, $0 \le t \le 1$, obviously shows that the induced map is surjective. The
hard part is the proof of injectivity which proceeds by showing that any loop is
homotopic to a concatenation of simple loops or even to a finite iteration of
a single simple loop.

For a long time it was an open problem whether higher homotopy groups of the
regular group $G(\c M)$ of a $W^{\a}$-algebra $\c M$ of type $II_1$ could 
serve in a finer classification of these algebras. However, this turned out not
to be the case as shown by the author in 1983 [Srd1] for complex and in 1985
[Srd2] for real $W^{\a}$-algebras. 

\medskip\lit{\bf Theorem 2}{\it If $\c M$ is a complex $W^{\a}$-algebra of
type $II_1$ then 
$$\pi_k(G(\c M)) \cong \cases
K_0(\Cal M) \,,\quad &k\equiv 1\, \pmod 2\cr
0\,,&k\equiv 0\, \pmod 2.\cr
\endcases$$
If $\c M$ is a purely real $W^\a$-algebra of type $II_1$ then
$$
\pi_k(G(\c M)) \cong \cases
K_0(\Cal M)\,,\quad &k\equiv 3\, \pmod 4\cr
0\,,&k\equiv 0,1,2\, \pmod 4.\cr
\endcases$$
}\par\medskip

\noindent
The proof relies essentially on a stabilization theorem similar to
Theorem 7
of section 3. It uses the fact that the regular group $G(\c M)$ of a type
$II_1$ $W^\a$-algebra $\c M$ is dense in $\c M$ which means that the
(topological) stable rank is equal to 1, and that any projection
$p$ in $\c M$ can be halved, i.e. $p = p_1 + p_2$ with two equivalent
orthogonal projections $p_1,p_2$. The periodicity comes from the
generalized Bott theorem. 

From these results one easily deduces the
homotopy type of homogeneous spaces re\-lative to a $W^\a$-algebra factor.
Conjectured by Breuer and proved in [EFT] the Grassmannian of a continuous
factor has simply connected components. The proof uses the fact mentioned
before, that any loop in $U(\c M)$ is homotopic to the iterate of a simple loop.
In [Srd3] we obtained the following conclusive results.
\medskip
\lit{\bf Theorem 3} {\it Let $\c M$
be a continuous $W^\a$-algebra factor, $Q = Q^\a = Q^{-1}\in \c M$, and, if $\c
M$ is purely real, $J^\a = - J = J^{-1} \in \c M$. By $U(\c M,Q) = U(Q,H_{\Bbb F})
\cap \c M$ we denote the group of $Q$-unitary elements relative to $\c M$ and by
$Sp(\c M,J) = Sp(J,H_{\Bbb R}) \cap \c M$ the symplectic group relative to $\c M$.
If, moreover, $p=1/2 (1+Q)$, $q = v^\a v$, where $J=v - v^\a$, then 
$$\pi_k
(U(\c M,Q)) \cong \pi_k(U(p\c Mp)) \oplus \pi_k(U((1-p)\c M(1-p)))\,,\quad k\ge
0$$ 
and 
$$ \pi_k (Sp(\c M,J)) \cong \pi_k(U(q\c Mq)_{\Bbb C})\,,\quad k\ge 0\,.$$
If $\c M$ is a type $II_1$-factor, $p$ a nontrivial projection then
the component ${Gr}_p(\c M)$ of the Grassmannian that contains $p$ has
homotopy groups
$$
\pi_k({Gr}_p(\c M)) \cong \cases
\Bbb R\,,\quad &k\equiv 0 \pmod{2d}\cr
0\,, &k\, \hbox{else},\cr
\endcases$$
where $d=1$ if $\c M$ is complex, and $d=2$ if $\c M$ is purely real.
The Banach homogeneous space $U_-(\c M)$ of skew-adjoint orthogonal elements in
a purely real  $II_1$-factor has homotopy groups
$$\pi_k(U_-(\c M)) \cong
\cases
\Bbb R\,,\quad &k\equiv 2\, \pmod 4\cr 0\,,&k\, \hbox{else}.\cr
\endcases$$}\par\medskip

Before we come back to the homotopy type of the regular group in other special
cases we want to take a look at $K$-theory. Most of the missing
references can be found in [Bla].
$K$-theory has entered operator theory in the mid-seventies when
Elliott used the dimension group ($= K_0(\c A)$) in order to classify $AF\,
C^\a$-algebras.
For example, the $UHF$-algebras $\c A_{(q)}$ defined as direct $C^\a$-limits
$(M(n_k,\Bbb F),i_k)$ of matrix algebras $M(n_k,\Bbb F)$, $\Bbb F = \Bbb R$ or $\Bbb C$ where
$n_{k+1} = q_kn_k$, $q_k \in \Bbb N$, and $i_k:M(n_k,\Bbb F) \ni x \mapsto x \otimes
1_{q_k} \in M(n_{k+1},\Bbb F)$, have $K$-groups
$$K_0(\c A_{(q)}) \cong \Bbb Z_{(q)} = \lim {\scriptstyle {1\over n_k}} \Bbb Z,\quad
K_1(\c A_{(q)}) = 0,$$ if $\Bbb F = \Bbb C$, and
$$KO_i(\c A_{(q)}) \cong \cases
\Bbb Z_{(q)}, \quad & i \equiv 0\, \pmod 4,\cr
\Bbb Z_2, \quad & i \equiv 1,2\, \pmod 8,\, 2^\infty \not| (q),\cr
0, \quad & i\, \hbox{else},\cr
\endcases$$
if $\Bbb F = \Bbb R$, using the fact [Han1] that the $K$-functor commutes with
direct $C^\a$-limits. Here $(q)$ denotes the formal product $\prod q_k$, and
$2^\infty \not|(q)$ means that 2 does occur at most finitely often in its formal
prime decomposition. The $AF\, C^\a$-algebras $\c A^{\theta}$ defined as direct
limits $(M(n_k,\Bbb F) \oplus M(m_k,\Bbb F),\phi_k)$, where $n_k$ and $m_k$ are defined
recursively by $n_1 = m_1 = 1$, $n_{k+1} = a_k n_k + m_k$, $m_{k+1} = n_k$
using the expansion of the irrational number $\theta$ as a continued fraction,
$\theta = [a_1;a_2,a_3,\dots ]$, and where $\phi_k(x,y) = (x \otimes 1_{a_k}
\oplus y,x)$, have $K$-groups $$K_0(\c A^{\theta}) \cong \Bbb Z + \theta \Bbb Z, \quad
K_1(\c A^{\theta}) = 0,$$ if $\Bbb F = \Bbb C$, and
$$ KO_i(\c A^{\theta}) \cong \cases
\Bbb Z + \theta \Bbb Z, \quad & i \equiv 0\, \pmod 4,\cr 
\Bbb Z_2^2, \quad & i \equiv 1,2\, \pmod 8,\cr
0, \quad & i\, \hbox{else},\cr
\endcases$$
if $\Bbb F = \Bbb R$ [Srd4].
The classification of real $AF\,C^\a$-algebras using real
$K$-theory has been achieved by Giordano [Gio].

$AF\,C^\a$-algebras are direct limits of zero-dimensional homogeneous
algebras. Direct limits $\c B_{(q)}$ of one-dimensional algebras
$(M(n_k,C(S^1,\Bbb F)),i_k)$ with $i_k$ a  $p_k$-times around imbedding have
been considered by Bunce and Deddens, cf. [Bla]. Analogously
to $UHF$-algebras one gets 
$$K_0(\c B_{(q)}) \cong \Bbb Z_{(q)}, \quad K_1(\c B_{(q)}) \cong \Bbb Z,$$ 
if $\Bbb F = \Bbb C$, and 
$$KO_i(\c B_{(q)}) \cong \cases
\Bbb Z_{(q)} (+ \Bbb Z_2), \quad
& i \equiv 0\, \pmod 8\, (\hbox{if } 2^\infty \not| (q)),\cr 
\Bbb Z_{(q)}, \quad & i \equiv 4\, \pmod 8,\cr
\Bbb Z_2 (+ \Bbb Z_2), \quad & i \equiv 1\, \pmod 8\, (\hbox{if } 2^\infty \not|
(q)),\cr 
\Bbb Z_2, \quad & i \equiv 2\, \pmod 8, \hbox{if } 2^\infty \not| (q),\cr
\Bbb Z, \quad & i \equiv 3\, \pmod 4,\cr
0, \quad & i\, \hbox{else},\cr
\endcases$$
if $\Bbb F = \Bbb R$ [Srd4]. For a large class of complex direct $C$*-limits of
one-dimensional algebras (those of real rank zero) Elliot has recently given 
also a classification using $K$-theory [Ell].

The first highlight in operator $K$-theory was obtained by Rieffel, Pimsner,
and Voiculescu who used the dimension group to distinguish (up
to Morita equivalence) between the irrational rotation algebras $\c A_{\theta}$
for various $\theta$, and, in particular, to exhibit nontrivial projections in
$\c A_{\theta}$ if $\theta \notin \Bbb Q$. The crucial device invented by Pimsner and
Voiculescu [PV1] for this purpose was the six-term exact sequence for crossed
products by $\Bbb Z$.

\medskip
\lit{\bf Theorem 4}{\it If $\c A$ is a $C^\a$-algebra provided with a
$\Bbb Z$-action $\alpha : \Bbb Z \to {Aut}(\c A)$ and $\c A \cp_{\alpha} \Bbb Z$ is the
corresponding crossed product then one has the cyclic exact sequence
$$
\matrix
K_0(\c A) &{\buildrel 1-\alpha_\a \over \to}&K_0(\c A)&{\buildrel 
i_\a \over \to}&K_0(\c A\times_\alpha \Bbb Z)\cr
\uparrow&&&&\downarrow\cr
K_1(\c A\cp_\alpha \Bbb Z)&{\buildrel i_\a \over \to}&K_1(\c A)&{\buildrel 
1-\alpha_\a \over \to}&K_1(\c A).\cr
\endmatrix$$ }\par\medskip

\noindent
Applying the Pimsner-Voiculescu sequence to the crossed product $\c A_{\theta}
= C(S^1,\Bbb C) \cp_{\alpha} \Bbb Z$, where $\Bbb Z$ acts on $C(S^1,\Bbb C)$ by rotating the
argument of a function by $e^{2\pi i \theta}$ gives the following
result.

\medskip\lit{\bf Corollary 1}{\it If $\theta \notin \Bbb Q$ and if $\c A_{\theta}
= C(S^1,\Bbb C) \cp_{\alpha} \Bbb Z$ is the corresponding irrational rotation algebra
then  $$K_0(\c A_{\theta}) \cong \Bbb Z + \theta \Bbb Z \quad \hbox{and} \quad K_1(\c
A_{\theta}) \cong \Bbb Z^2.$$}\par\medskip

\noindent
An analogous cyclic exact sequence with 24 terms holds in the real case [Ros3],
[Srd4], and [Sty] and gives:

\medskip\lit{\bf Corollary 2}{\it If $\theta \notin \Bbb Q$ and if $\c
A_{\theta}^{\Bbb R} = C(S^1,\Bbb R) \cp_{\alpha} \Bbb Z$ is the corresponding real irrational
rotation algebra then
$$KO_i(\c A_{\theta}^{\Bbb R}) \cong \cases
\Bbb Z + \theta\Bbb Z + \Bbb Z_2,  & i \equiv 0\, \pmod 8,
\cr \Bbb Z + \Bbb Z_2^3, & i \equiv 1\, \pmod 8, \cr
\Bbb Z_2^3, & i \equiv 2\, \pmod 8, \cr
\Bbb Z + \Bbb Z_2, & i \equiv 3\, \pmod 8, \cr
\Bbb Z + \theta\Bbb Z, & i \equiv 4\, \pmod 8, \cr
\Bbb Z, & i \equiv 5,7\, \pmod 8,
\cr0, & i \equiv 6\, \pmod 8. \cr
\endcases$$}\par\medskip

\noindent
To obtain these corollaries one has to notice that the maps $\alpha_\a$ induced
by the action are the identity (because the phase $e^{2\pi i \theta}$ can be
deformed to 1 and the $K$-groups are not affected by the deformation [AP]) 
and to prove that the so obtained short exact
sequences split. Note that the isomorphism of $K_0$ with a subgroup
of $\Bbb R$ is induced by the unique tracial state of these algebras; cf. [PV1],
[Pim1]. Crossed products coming from torus actions are treated in [Rou].

Just as $\c A_\theta$ can be considered as a deformation of the
commutative $C^\ast$-algebra $C(T^2)$ one can deform the relations that 
determine, e.g., $S^3 = SU(2)$ as a submanifold of $\Bbb C^2$ by inserting a 
phase factor.
Thus, for $\theta \notin \Bbb Q$ the noncommutative 3-sphere $S^3_\theta$ can
either be defined by suitable generators and relations or more
geometrically by gluing two noncommutative solid tori [Mat]:
For $D = \{z \in \Bbb C \mid |z|\le 1\}$ let $\c D_\theta = C(D)
\cp_\alpha \Bbb Z$ with $\alpha$ as in Corollary 1. Then one has a canonical
surjection $\pi_\theta : \c D_\theta \to \c A_\theta$. Define an
isomorphism $\rho : \c A_{-\theta} \to \c A_\theta$ by $\rho(\hat u) =
v$ and $\rho(\hat v) = u$ on the corresponding generators and let
$$S^3_\theta = \{ (a,b) \in \c D_{-\theta}\oplus \c D_\theta \mid
\rho(\pi_{-\theta}(a)) =\pi_\theta(b)\}.$$
The computation of the corresponding $K$-groups uses the following Mayer-Vietoris
sequence for a pull-back
$$\matrix 
\c D & {\buildrel {g_1} \over \longrightarrow} &
\c B_1 \cr
\llap{g_2}\downarrow&&\downarrow\rlap{f_1}\cr
\c B_2 & {\buildrel {f_2} \over \longrightarrow}  & \c C \cr
\endmatrix$$
of $C^\ast$-algebras
where $\c D =\{ (b_1,b_2) \in \c B_1 \oplus \c
B_2 \mid f_1(b_1) = f_2(b_2)\}$, $f_1$, $f_2$ are surjective, and where
$g_i :\c D \to \c B_i$, $i=1,2$, are
the restrictions of the canonical projections, cf. [Bla] and
[Srd4] (also for more applications).
\medskip
\lit{\bf Theorem 5}{\it Any pull-back diagram of $C^\a$-algebras gives rise
to a long exact sequence
$$\to K_n(\c D) {\buildrel {(g_{1 \a},g_{2 \a})} \over \longrightarrow}
 K_n(\c B_1) \oplus K_n(\c B_2) {\buildrel {f_{2 \a}-f_{1 \a}} \over
\longrightarrow} K_n(\c C) {\buildrel {\alpha} \over \longrightarrow}
K_{n-1}(\c D) \to .$$}\par\medskip

\medskip
\lit{\bf Corollary 1}{\it If $\theta \notin \Bbb Q$ then the K-groups of
$S^3_\theta$ are
$$K_0(S^3_\theta) \cong \Bbb Z \cong K_1(S^3_\theta).$$
}\par\medskip

The real case is more involved. Here one has to ``glue'' real subalgebras of 
irrational rotation algebras with different real structure. The real structure 
is defined on the generators by $\bar u=u^\a$, $\bar v=v$ in case of 
$\c A_\theta$ and by $\bar{\hat v}={\hat v}^\a$ and $\bar{\hat u}=\hat u$ 
in case of $\c A_{-\theta}$. Now, since the latter $C^\a$-algebra is not a 
crossed-product, the corresponding real form for $\c D_{-\theta}$ is the 
``cone'' $\{f\in C([0,1],\c A_{-\theta}^{\Bbb R}\mid f(0)\in C(S^1,\Bbb R)\}$. Then 
one obtains for the real noncommutative 3-sphere $S_\theta^{3\Bbb R}$:
\medskip
\lit{\bf Corollary 2}{\it If $\theta \notin \Bbb Q$ then the K-groups of
$S_\theta^{3\Bbb R}$ are
$$K_0(S_\theta^{3\Bbb R}) \cong \cases
\Bbb Z + \Bbb Z_2,  & i \equiv 0\, \pmod 8,
\cr \Bbb Z_2^2, & i \equiv 1\, \pmod 8, \cr
\Bbb Z_2, & i \equiv 2\, \pmod 8, \cr
\Bbb Z, & i \equiv 3,4,~\hbox{or}~7\, \pmod 8, \cr
0, & i \, \hbox{else}. \cr
\endcases$$}\par\medskip

More generally, one can define noncommutative lens spaces $L_\theta(p,q)$ using 
the gluing isomorphism that is given by 
$$\rho_A(u)=e^{-k\ell\pi i\theta}v^k u^\ell=\hat 
v\quad\hbox{and}\quad\rho_A(v)=e^{-pq\pi i\theta}v^qu^p=\hat u$$ with 
$A=\pmatrix
q&k\cr 
p&\ell\cr
\endpmatrix\in SL(2,\Bbb Z)$. Slightly more elaborate calculations yield
$$K_0(L_\theta(p,q)) \cong \Bbb Z_p+\Bbb Z, \quad K_1(L_\theta(p,q) \cong 
\Bbb Z,$$                           
$$K_0(L_\theta(0,1)) \cong \Bbb Z^2, \quad K_1(L_\theta(p,q) \cong 
\Bbb Z^2,$$ 
in the complex case (cf.\ [MT]), and 
$$KO_i(L_\theta^{\Bbb R}(p,q)) \cong \cases
\Bbb Z + \Bbb Z_p + \Bbb Z_2, \quad & i \equiv 0\, \pmod 8 \cr 
\Bbb Z_2^2 + \Bbb Z_{(p,2)}, \quad & i \equiv 1\, \pmod 8,\cr
\Bbb Z_2 + \Bbb Z^2_{(p,2)}, \quad & i \equiv 2\, \pmod 8, \cr 
\Bbb Z + \Bbb Z_{(p,2)}, \quad & i \equiv 3\, \pmod 8,\cr
\Bbb Z + \Bbb Z_p, \quad & i \equiv 4\, \pmod 8,\cr
\Bbb Z, \quad & i \equiv 7\, \pmod 8,\cr
0, \quad & i\, \hbox{else},\cr
\endcases$$
in the real case [Sdr5].

As a final application of the Mayer-Vietoris sequence we give the $K$-groups 
of $C^\a$-algebras associated with compact orbifolds of dimension 2, i.e., 
singular algebraic curves. Topologically such a curve is an oriented compact 
surface of genus $g$ whose singular locus consist of a finite set of points of 
multiplicity $n_j$, i.e.\ at such a point $n_j$ branches meet. The complex case 
has been treated by Carla Farsi [Far] the real case has been determined by the 
author [Srd5]. Denoting by $\ell$ the number of 
singular points, and $\sigma=1+\sum_{j=1}^\ell(n_j-1)$ one obtains
\medskip
\lit{\bf Corollary 3}{\it For an orbifold $T_{g,\ell}$ of genus $g$ 
with $\ell$ singular points of multiplicity $n_j$ the K-groups of the 
corresponding $C^\a$-algebras $C^\a(T_{g,\ell})$ and $C^\a(T_{g,\ell}^{\Bbb 
R})$ are given by
$$K_0(C^\a(T_{g,\ell})) \cong \Bbb Z^{1+\sigma}\quad\hbox{and}\quad 
K_1(C^\a(T_{g,\ell}))\cong \Bbb Z^{2g}$$
and by
$$KO_i(C^\a(T_{g,\ell}^{\Bbb R})) \cong \cases
\Bbb Z^\sigma + \Bbb Z_2^{2g+1}, \quad & i \equiv 0\, \pmod 8 \cr 
\Bbb Z_2^{2g+\sigma}, \quad & i \equiv 1\, \pmod 8,\cr
\Bbb Z + \Bbb Z_2^\sigma, \quad & i \equiv 2\, \pmod 8, \cr 
\Bbb Z^{2g}, \quad & i \equiv 3\, \pmod 8,\cr
\Bbb Z^\sigma, \quad & i \equiv 4\, \pmod 8,\cr
0, \quad & i \equiv 5\, \pmod 8,\cr
\Bbb Z, \quad & i \equiv 6\, \pmod 8,\cr
\Bbb Z^{2g} + \Bbb Z_2, \quad & i \equiv 7\, \pmod 8,\cr
\endcases$$
}\par\medskip

Another class of examples is obtained by extensions, i.e.
$C^\ast$-algebras $\c A$ that fit into a short exact sequence
$$0 \to \c J \to \c A \to \c B \to 0.$$
From this sequence one obtains a long exact sequence of $K$-groups (of
period 6 in the complex and of period 24 in the real case)
$$\to K_n(\c J) \to K_n(\c A) \to K_n(\c B) \to K_{n-1}(\c J) \to\, .$$
As an application we consider the Toeplitz $C^\ast$-algebra $\c
T(S^{2n-1})$ and the Calderon-Zygmund algebra $CZ(S^k)$ of zero-order
pseudodifferential operators on $S^k$ which are both
extensions of the compacts by a commutative
$C^\a$-algebra, viz.
$$0 \to \c K \to \c T(S^{2n-1}) \to C(S^{2n-1}) \to 0,$$
and
$$0 \to \c K \to CZ(S^{k}) \to C(S^{\a}S^{k}) \to 0,$$
where $S^{\a}S^{k}$ is the cosphere bundle over $S^{k}$.
Now $K_n(C(S^k)) = K_n(\Bbb C) + K_{n+k}(\Bbb C)$ and
$$K_0(C(S^\ast S^k)) \cong \cases
\Bbb Z + \Bbb Z_2, \, &k\, \hbox{even}\cr
\Bbb Z + \Bbb Z,\, &k\, \hbox{odd}\cr
\endcases
\quad\hbox{and}\quad
K_1(C(S^\ast S^k)) \cong \cases
\Bbb Z, \, &k\, \hbox{even}\cr
\Bbb Z + \Bbb Z,\, &k\, \hbox{odd},\cr
\endcases$$
so that
$$K_0(\c T(S^{2n-1})) \cong \Bbb Z \cong K_1(\c T(S^{2n-1}))$$
and [Les]
$$\eqalign{K_0(CZ(S^{2n})) &\cong \Bbb Z + \Bbb Z_2,\quad K_1(CZ(S^{2n})) = 0,\cr
K_0(CZ(S^{2n-1})) &\cong \Bbb Z^2,\quad K_1(CZ(S^{2n-1})) \cong \Bbb Z.\cr}$$ 
These groups can in turn be used to compute the $K$-groups of other
Toeplitz $C^\ast$-algebras.
For the $C^\a$-algebra $\c T(S^{2n_1-1} \times \cdots \times S^{2n_k-1})$ of
Toeplitz operators on the Hardy space $H^2(S^{2n_1-1} \times \cdots \times 
S^{2n_k-1})$ with continuous symbols or $\c T(L_n)$ the corresponding algebra
on the Hardy space over the Lie spheres $L_n$ the $K$-group $K_0(\c T)$ is
isomorphic to $\Bbb Z$ and $K_1(\c T)$ is trivial [Deu1], [Les]. In the
second case one has again an extension
$$0 \to \c K(H^2(S^1))\otimes CZ(S^{n-1}) \to \c T(L_n) \to C(L_n)\to
0$$
and uses
$$K_0(C(L_k)) \cong \cases
\Bbb Z + \Bbb Z, \, &k\, \hbox{even}\cr
\Bbb Z,\, &k\, \hbox{odd}\cr
\endcases\quad\hbox{and}\quad
K_1(C(L_k)) \cong \cases
\Bbb Z + \Bbb Z, \, &k\, \hbox{even}\cr
\Bbb Z + \Bbb Z_2,\, &k\, \hbox{odd},\cr
\endcases$$
while the first is a direct consequence of the K\"unneth short split
exact sequence
$$ 0 \to K_\ast(\c A)\otimes K_\ast(\c B) \to K_\ast(\c A\otimes \c B)
\to Tor_1^{\Bbb Z}(K_\ast(\c A),K_\ast(\c B))\to 0.$$
 
The proof given in [Srd4] and in [Sty] for the real Pimsner-Voiculescu sequence and given by
Joachim Cuntz in [Cun2,5] as an alternative proof in the complex case is based
on the homotopy properties of the regular group of another class of
$C^{\a}$-algebras, the so-called Cuntz algebras $\c O_n$, $n \ge 2$,
first defined by Cuntz in 1977. 

\medskip\lit{\bf Theorem 6}{\it Let $\c O_n^{(\Bbb R)}$ be the (real)
$C^{\a}$-algebra generated by n (real) partial isometries $S_j$ with
$S_j^{\a}S_j = 1$ and $\sum\limits_{j=1}^n S_j S_j^\a = 1$, then 
$$
K_i(\c O_n) \cong \cases
\Bbb Z_{n-1}\,,\quad &i\equiv 0\, \pmod 2,\cr
0\,,\quad &i\equiv 1\, \pmod 2,\cr
\endcases$$
and
$$
KO_i(\c O^{\Bbb R}_n) \cong \cases
\Bbb Z_{n-1}\,,\quad &i\equiv 0\, \pmod 4\cr
\Bbb Z_2\,,\quad &i\equiv 1,3, \pmod 8, n\, \hbox{odd},\cr
\Bbb Z_2^2\,,\quad &i\equiv 2\, \pmod 8, n\, \hbox{odd},\cr
0\,,\quad &i\, \hbox{else}\cr
\endcases$$
}\par\medskip

\noindent
We refer to [Cun3] and [Srd4] for the $K$-groups of the more general 
Cuntz-Krieger $C^\ast$-algebras.

A generalization of Bott's periodicity theorem has been obtained by
Connes for actions $\alpha : \Bbb R \to {Aut}(\c A)$ [Con] - see [Ksp2] and
[Srd4] for the proof of the real version.
\medskip\lit{\bf Theorem 7}{\it For
any complex $C^{\a}$-algebra $\c A$ and any $\Bbb R$-action $\alpha$ on $\c A$ one
has isomorphisms $$
K_i(\c A\cp_\alpha \Bbb R) \cong K_{1-i} (\c A)\,,\quad i=0,1\,.$$
If $\c A$ is a real $C^{\a}$-algebra there are isomorphisms
$$
KO_i(\c A\cp_\alpha \Bbb R) \cong KO_{i-1}(\c A)\,,\quad i\in \Bbb Z\,.$$}\par\medskip

\noindent
His main purpose was to compute $K$-groups of group $C^{\a}$-algebras $C^\a(G)$
and $C^\a_{red}(G)$ where $G$ is a connected Lie group and $C^\a(G)$ 
denotes the completion of the convolution algebra
$L^1(G)$ with respect to the greatest $C^\a$-norm and $C^\a_{red}(G)$
with respect to the norm induced by
left regular representation on $L^2(G)$. The real group $C^\a$-algebras
$C^\a(G,\Bbb R)$ and $C^\a_{red}(G,\Bbb R)$ are the closure of $L^1(G,\Bbb R)$ under the
respective norm in $C^\a(G)$ and $C^\a_{red}(G)$. For a compact group $G$
one has  $$K_0(C^{\a}(G)) \cong \Bbb Z \hat G = R(G), \quad K_1(C^\a(G)) = 0,$$
$R(G)$ the representation ring of $G$, and, more generally, for an action
$\alpha$ of $G$ on a $C^\a$-algebra $\c A$ the so-called Green-Julg theorem
holds: $$K_i(\c A \cp_{\alpha}G) \cong K_i^G(\c A), \quad i = 0,1.$$
For $G$ a noncompact connected Lie group with maximal compact subgroup $H$ 
Connes conjectured isomorphisms [BaC]
 $$K_i(C_{red}^\a (G)) \to K_H^{i+\dim G/H}(\{pt\}), \quad i = 0,1.$$ 
This conjecture has since been confirmed in several special cases. By
successive application of Theorem 7 one gets [Con] for a connected
simply connected solvable Lie group $G$ that 
$$
K_i(C_{red}^\a (G)) \cong K^{i+\dim G}(\{pt\}) \cong \cases
\Bbb Z,\quad &
i \equiv\, {\dim} G \pmod 2,\cr
0,\quad &i \not\equiv\, {\dim} G \pmod 2,\cr
\endcases$$
or more generally for a $G$-$C^\a$-algebra $\c A$
$$K_i(\c A \cp G) \cong K_{i+\dim G}(\c A).$$

\noindent
The conjecture is also true 
if $G$ is nilpotent, or a motion group, i.e.\ $G = V \cp K$, $V$ a
vector group and $K$ a connected compact Lie group acting linearly on $V$, or
if $G$ is semi-simple of rank 1 [Ros1,2].

\noindent
If $G$ is a connected reductive linear group and $H$ a
maximal compact subgroup then 
$$K_i(C_{red}^\a (G)) \cong \cases
R_{spin}(H),\quad &i \equiv 
{\dim} (G/H)\, \pmod 2,\cr 0,\quad &i\, \hbox{else} ,\cr
\endcases$$
where $R_{spin}(H)$ is the $R(H)$-module of ``spinoral'' representations of $H$.
This has been proved in special cases by Penington and Plymen
and by Valette, and in full generality by Wasserman
[Was].

Another important class of $C^\a$-algebras are the group $C^\a$-algebras of
discrete groups. Again a great impact has been made by Pimsner and Voiculescu
who considered crossed products by actions of finitely generated free groups.
There is a series of papers devoted to the $K$-theory of such group
$C^\a$-algebras and of their crossed pro\-ducts by Anderson and Paschke,
Cuntz, Lance, Kasparov, Natsume, and finally by
Pimsner [Pim2]. Pimsner's
paper contains all the previous ones as special cases. The outcome is a kind of
Mayer-Vietoris sequence for the $K$-groups of the group $C^\a$-algebra of a group
acting without involution on a tree and whose building blocks are the
stabilizers of edges and vertices of the tree. We do not state the exact
sequence here but emphasize that it holds in the real and in the complex case,
[Srd5]. Using this exact sequence the $K$-groups of many group $C^\a$-algebras
can be computed.
$$
\matrix
K_i(C^\ast_{red}(\Gamma))\,\Gamma &  F_k & \Bbb Z_k\ast _{\Bbb Z_{\ell}}\Bbb Z_m &
\Gamma_g & \Sigma_k & H_d & \Bbb Q\ast \Bbb Q  \cr
i\cr
0 &  \Bbb Z &  \Bbb Z^{k+m-\ell} &  \Bbb Z^2 & \Bbb Z &  \Bbb Z^3 &  \Bbb Z  
 \cr   
1 &  \Bbb Z^k &  0 &  \Bbb Z^{2g} &  \Bbb Z^{k-1}+\Bbb Z_2 &
\Bbb Z^3 &  \Bbb Q^2 \cr
KO_i(C^\ast_{red}(\Gamma))\cr
i\cr  
0 &  \Bbb Z &  \Bbb Z^{k+m-\ell} &  \Bbb Z &  \Bbb Z & \Bbb Z &  \Bbb Z 
\cr   
1 &  \Bbb Z^k + \Bbb Z_2 &  \Bbb Z_2^{k+m-\ell} &  \Bbb Z^{2g}+\Bbb Z_2
&  \Bbb Z^{k-1}+\Bbb Z_2^2 &  \Bbb Z^2+\Bbb Z_2 &  \Bbb Q^2+\Bbb Z_2 \cr  
2 &  \Bbb Z_2^{k+1} &  \Bbb Z_2^{k+m-\ell} &  \Bbb Z+\Bbb Z_2^{2g+1} & 
\Bbb Z_2^{k+1} & \Bbb Z^2+\Bbb Z_2^3 &  \Bbb Z_2 \cr  
3 &  \Bbb Z_2^k &  0 &  \Bbb Z_2^{2g+1} &  \Bbb Z_2^{k+1} &  \Bbb Z+\Bbb Z_2^4 & 
0\cr  
4 &  \Bbb Z & \Bbb Z^{k+m-\ell} &  \Bbb Z+\Bbb Z_2 &  \Bbb Z+\Bbb Z_2 &  \Bbb Z+\Bbb Z_2^3 & 
\Bbb Z \cr  
5 &  \Bbb Z^k &  0 &  \Bbb Z^{2g} &  \Bbb Z^{k-1}+\Bbb Z_2 &  \Bbb Z^2+\Bbb Z_2 & 
\Bbb Q^2 \cr  
6 &  0 &  0 &  \Bbb Z &  0 &  \Bbb Z^2 & 0 \cr  
7 &  0 &  0 &  0 & 0 &  \Bbb Z &  0 \cr 
\endmatrix$$ 
Here $F_k$ denotes the free group on $k$ generators, $\Gamma_g$ and $\Sigma_k$
are the fundamental groups of oriented resp.\ nonorientable compact surfaces,
and $H_d$ is the discrete Heisenberg group defined as before but with $x,y,z
\in \Bbb Z$. We have only listed $K$-groups of the reduced $C^\a$-algebras. But in
all cited examples they coincide with the $K$-groups of the corresponding full
$C^\a$-algebras, because the groups are $K$-amenable, a notion
introduced by Cuntz for discrete
groups and by Julg and Valette for locally compact
continuous groups. Using different methods Kasparov and Skandalis [KSk]
obtain a still more general result for groups operating on a
Bruhat-Tits building.

We also like to mention that $K$-theory helped to confirm an old standing
conjecture by Kadison on the existence of nontrivial idempotents in
the group $C^\a$-algebra of a torsion-free discrete group in many special
cases. This conjecture is related to another one: If $\Gamma$
is torsion-free, then the canonical trace $\tau$ induces a surjection $\tau_\a
: K_0(C^\a_{red}(\Gamma)) \to \Bbb Z$. While this map is in general not
one-to-one, e.g.\ for $\Gamma = \Gamma_g$, our calculations indicate the
following strengthening: If $\Gamma$
is torsion-free, then the canonical trace $\tau$ restricted to the real
group $C^\a$-algebra induces an isomorphism $\tau_\a
: KO_0(C^\a_{red}(\Gamma,\Bbb R)) \to \Bbb Z$.

Now we come back to nonstable $K$-theory, i.e.\ the homotopy groups of the
regular group of a $C^\a$-algebra $\c A$ itself. In the previous section we have
already noted a connection between the size of the stable range
$$SR(\c A) = \{(n,k) \in \Bbb Z_+ \times \Bbb N |\, k \ge k(n)\},$$
with $k(n) = \min \{k | \pi_n(GL(k,\c A)) \cong \pi_n(GL(\infty,\c A))\}$, and
 the stable rank $sr(\c A)$ if this is finite.
Elaborating on the proof of Theorem 2, Rieffel
found the following general results [Rie2].
\medskip
\lit{\bf Theorem 8} {\it If $\c A$ is a unital $C^\a$-algebra, $p \in \Bbb N$, and
${csr}(C(T^k,\c A)) \le p$ for all $k \ge 0$, then
$$\pi_k(GL(n,\c A)) \cong \cases
K_1(\Cal A), \quad  & k\, \hbox{even}, \cr K_0(\Cal A),
\quad & k\, \hbox{odd},\cr
\endcases$$
for all $n \ge p-1$ and $k \ge 0$.}\par\medskip
\noindent 
Here ${csr}(\c A)$ denotes the connected stable rank,
i.e.\ the least integer $m$ such that $GL^0(n,\c A)$ acts transitively on 
$$Lg(n,\c A) = \{(a_1,\dots,a_n)\in \c A^n | \sum b_ia_i = 1 \hbox{ for some }
(b_1,\dots,b_n)\in \c A^n\}$$
for all $n \ge m$.
 
\noindent
As corollaries he obtains that $\c A$ has full stable range
if $\c A$ is ${tsr}$-boundedly divisible, i.e. if there is a constant $c$,
such that for every $m$ there is an $n \ge m$ and a $C^\a$-algebra $\c B$ with
$\c A \cong M(n,\c B)$ and ${tsr}(\c B) \le c$.
Examples  of trs-boundedly divisible $C^\a$-algebras are type $II_1$ 
$AW^\a$-algebras and tensor products $\c A \otimes \c B$ where $\c B$ 
is a unital divisible $AF\,C^\a$-algebra (e.g. a $UHF$-algebra) and $\c A$ 
is unital with ${tsr}(\c A) < \infty$ or 
$\c A = C(X)\cp_\alpha \Gamma$ is a crossed product with a discrete solvable 
group $\Gamma$ acting on a compact space $X$. Thomsen has shown [Tho] that 
$\c A \otimes \c B$
has also full stable range if $\c A$ is a unital $C^\a$-algebra and $\c B$ is 
an $AF$ $C^\a$-algebra with $K_0(\c B)$ having large denominators in the 
sense of
[Nis], in particular, if $\c B$ is infinite dimensional and simple.

\noindent
Moreover, he considers the groups of quasi-invertibles and of quasi-unitaries,
$G_q(\c A)$ resp. $U_q(\c A)$, in a nonunital complex $C^\a$-algebra $\c A$, and
shows that an exact sequence
$$0 \to \c J {\buildrel i \over \to} \c A {\buildrel p \over \to} \c B \to 0 $$
of $C^\a$-algebras gives rise to a long exact sequence of homotopy groups
$$\to \pi_k(U_q(\c J)) {\buildrel i_{\a} \over \to} \pi_k(U_q(\c A))
{\buildrel p_{\a} \over \to} \pi_k(U_q(\c B)) \to \pi_{k-1}(U_q(\c
J))\to\cdots\to \pi_0(U_q(\c B)).$$
Note that an analogous exact sequence can be obtained for the groups of
quasi-invertibles. Then the proof uses the selection theorem of Michael
at a crucial point, and leads to the same conclusion for an exact
sequence of real $C^\a$-algebras. If $X$ is a compact space Thomsen obtains the
following partial result ($C(X) = C(X,\Bbb C)$). \medskip 
\lit{\bf Theorem 9}{\it
If $X$ is a compact Hausdorff space then 
$$\pi_n(U(C(X))) \cong H^{1-n}(X,\Bbb Z),\quad n \ge 0.$$ 
\lit{}In particular, $\pi_1(U(C(X))) \cong \Bbb Z,$ if $X$ is connected, and
$\pi_n(U(C(X))) = 0$ for $n \ge 2$.}
\par\medskip

\noindent
As an immediate corollary we obtain the homotopy groups of the regular group
of the Toeplitz $C^\a$-algebra $\c T(S^{2n-1})$ and of the Calderon-Zygmund
algebra $CZ(S^{n})$. \medskip
\lit{\bf Corollary}{\it For $n \ge 1$ and $\c A = \c T(S^{2n-1})$ or 
$ CZ(S^{n+1})$ one has 
$$\pi_k(G(\c A)) = \cases
\Bbb Z, \quad &k \equiv 1\, \pmod 2, k \ne 1,\cr
0, \quad &k \equiv 0\, \pmod 2,\cr
\Bbb Z^2, \quad & k = 1,\cr
\endcases$$
while
$$\pi_k(G(CZ(S^{1}))) = \cases
\Bbb Z, \quad &k \equiv 1\, \pmod 2, k \ne 1\, \hbox{or } k
= 0,\cr 0, \quad &k \equiv 0\, \pmod 2, k \ne 0,\cr
\Bbb Z^3, \quad & k = 1.\cr
\endcases$$}\par\medskip

\noindent
The results for $CZ(S^1)$ can also be found in the papers by
Khimshiashvili [Khi1,2]. Now in both cases the stable rank is equal to $n-1$ or
to $n$, and there are homotopy groups of the regular group of matrix algebras
which have torsion. To determine these homotopy groups is as hard as the
computation of the homotopy groups of classical groups in the instable range.

Now the relation between stable rank and size of the stable range of a
$C^\a$-algebra is not completely settled. There are many examples of extremely noncommutative
$C^\a$-algebras with infinite stable rank and full stable range, e.g. purely 
infinite simple
$C^\a$-algebras (such as the Cuntz algebras $\c O_A$) [Zha1-4] or $\c A
\otimes \c B$, where $\c B$ is an arbitrary $C^\a$-algebra (in the nonunital case
consider quasi-invertibles), and $\c A$ is either $\c O_n^{(\Bbb R)}$ [Tho]
([Srd4]), a $W^\a$-algebra $\c M$ of infinite type, the multiplier algebra
$M(\c D \otimes \c K)$ of $\c D \otimes \c K$, where $\c D$ is a
 $C^\a$-algebra with a strictly positive element, or the corresponding
Calkin algebra [Tho].

So it finally appears that the
crucial property for full stable range is some kind of divisibility which is
guaranteed if the algebra has a lot of projections. And if there are enough 
projections not
equivalent to the identity, the regular group tends to be contractible, cf.\
section 2. The most important example here is the algebra 
$M(\c A \otimes \c K) \otimes \c B$ where $\c A$ has a strictly positive
element and $\c B$ is arbitrary. This generalization of Kuiper's theorem (which
is $\c A = \Bbb C = \c B$) has been proved by Mingo, Kasimov, and Troitskii
in special cases, by Cuntz and Higson [CH] when $\c B = \Bbb C$, and by Thomsen
[Tho] and Troitskii [Tro] in general (see [Srd4] for real $C^\a$-algebras).

As in
the case $\c A = \Bbb F$ the general Kuiper theorem has applications to classifying
spaces using generalized Fredholm operators, see [Min] and [Srd4]. Fredholm
operators relative to the ideal of ``compact'' operators in a semifinite
$W^\a$-algebra (or $AW^\a$-algebra) have been considered by [Bre], [Gew],
[Php], [Rae], and [Srb], see also [CPh].

\vfill\break

\parindent=13pt

\centerline{\bf References}
\bigskip\bigskip

\ref{AC} {Alling, N.L., Campbell, L.A.} {Real Banach algebras II}{Math.\
Z.\ 125 (1972) 79 - 100}

\ref{AP} {Anderson, J., Paschke, W.} {The rotation algebra}{Houston Math.\ J.\ 
15 (1989) 1-26}

\ref{ASS} {Araki, H., Smith,
M.-S.B., Smith, L.} {On the homotopical significance of the type of von Neumann
algebra factors}{Comm.\ Math.\ Phys.\ 22 (1971) 71 - 88}

\ref{Ati}{Atiyah, M.F.}{K-Theory}{W.A. Benjamin, New York, 1967}

\ref{AS}  {Atiyah, M.F., Singer, I.M.} {Index theory for skew-adjoint
Fredholm operators}{ Publ.\ Math.\ IHES 37 (1969) 5 - 26 }

\ref{BaC}{Baum, P., Connes, A.}{Geometric K-theory for Lie groups and
foliations}{Preprint, IHES, 1982}

\ref{Bel} {Belov, I.S.} {The homotopy type of the linear group of the
Banach space $C(\Gamma_{m\omega_1})$}{Amer.\ Math.\ Soc.\ Transl.\ (2) 115
(1980) 169 - 174}

\ref{Bla}   {Blackadar, B.E.} {K-Theory for Operator Algebras}
{Springer-Verlag, New\- York - Berlin, 1986}

\ref{Blu} {Blum, E.K.} {The fundamental group of the principal component of
a commutative Banach algebra}{Proc.\ Amer.\ Math.\ Soc.\ 4 (1953) 397 -
400}

\ref{BH}   {Borel, A., Hirzebruch, F.} {Characteristic classes and
homogeneous spaces II}{Amer.\ J.\ Math.\ 81 (1959) 315 - 382 }

\ref{Bos}{Bost, J.-B.}{Principe d'Oka, K-th\'eorie et syst\`emes dynamiques
non commutatifs}{Invent.\ Math.\ 101 (1990) 261 - 333}

\ref{Bot}   {Bott, R.H.} {The stable homotopy of the classical groups}
{Proc.\ Natl.\ Acad.\ Sci.\ USA 43 (1957) 933 - 935, Ann.\ of Math.\ 70 (1959)
313 - 337}

\ref{Bre}  {Breuer, M.} {Theory of Fredholm operators and vector bundles
relative to a von Neumann algebra}{Rocky Mountain J.\ Math.\ 3 (1973) 383 - 429
}

\ref{BW}  {Br\"uning, J., Willgerodt, W.} {Eine
Verallgemeinerung eines Satzes von N. Kui\-per}{Math.\ Ann.\ 220 (1976) 47 - 58
}

\ref{CE} {Carey, A.L., Evans, D.E.} {Algebras almost commuting with Clifford
algebras} {J. Funct. Anal. 88 (1990) 279 - 298}

\ref{CPh}{Carey, A.L., Phillips, J.}{Algebras almost commuting with Clifford
algebras in a $II_\infty$ factor}{K-Theory 4 (1991) 445 - 478}

\ref{Ctn}   {Cartan, E.} {Sur les nombres de Betti des espaces de groupes
clos}{C.\ R.\ Acad.\ Sci.\ 187 (1928) 196-198, see also l'Enseignement
Math.\ 35 (1936) 177 - 200 }

\ref{Con}  {Connes, A.} {An analogue of the Thom isomorphism for
crossed products of a $C^\a$-algebra by an action of $\Bbb R$}{Adv.\ in Math.\ 39
(1981) 31 - 55}

\ref{Cor} {Corach, G.} {Homotopy stability 
in Banach algebras}{Rev.\ Union Mat.\ Argent.\ 32 (1986) 233 - 243}

\ref{CoL}  {Corach, G., Larotonda, A.R.} {A stabilization theorem for
Banach algebras}{ J.\ Algebra 101 (1986) 433 - 449}

\ref{CPR}{Corach, G., Porta, H., Recht, L.}{Differential geometry of systems of
projections in Banach algebras}{Pac.\ J.\ Math.\ 143 (1990) 209 - 228}

\ref {CoS1}{Corach, G., Su\'arez, F.D.}{Dense morphisms in commutative Banach
algebras}{Trans.\ Amer.\ Math.\ Soc.\ 304 (1987) 537 - 547}

\ref {CoS2}{Corach, G., Su\'arez, F.D.}{Continuous selections and stable rank
of Banach algebras}{Topology Appl. 43 (1992) 237 - 248}

\ref {CrL}{Cordes, H.O., Labrousse, J.P.}{The invariance of the index in the
metric space of closed operators}{J.\ Math.\ Mech.\ 12 (1963) 693 - 720}

\ref {Cue}{Cuellar, J.R.}{Fredholmoperatoren auf
lokal beschr\"ankten R\"aumen mit Anwendungen auf elliptische
Gleichungen}{Dissertation, Univ.\ Mainz, 1982}

\ref{Cun1}  {Cuntz, J.} {K-theory for certain $C^\a$-algebras}{Ann.\ of
Math. 113 (1981) 181 - 197}

\ref{Cun2}  {Cuntz, J.} {K-theory for certain $C^\a$-algebras II}{J. Operator
Theory 5 (1981) 101 - 108}

\ref{Cun3}{Cuntz, J.}{A class of $C^\a$-algebras and topological Markov chains
II}{Invent.\ Math.\ 63 (1981) 25 - 46}

\ref{Cun4} {Cuntz, J.} {K-theoretic amenability for discrete groups}{J.\
reine angew.\ Math.\ 344 (1983) 180 - 195}

\ref {Cun5}{Cuntz, J.}{$K$-theory and $C^\a$-algebras}{In: Algebraic K-Theory,
Number Theory, Geometry and Analysis, Bielefeld, 1982, Lect.\
Notes Math.\ 1046, pp.\ 55 - 79, Springer-Verlag, Berlin, 1984}
 
\ref{CH}   {Cuntz, J., Higson, N.} {Kuiper's theorem for Hilbert modules}{In:
Operator Algebras and Mathematical Physics, Iowa City, 1985, Con\-temp.\ Math.\
62, pp.\ 429 - 435, Amer.\ Math.\ Soc., Providence, R.I., 1987}

\ref{vDa}  {Daele, A. van} {K-theory for graded Banach algebras I,II}{
Quart.\ J.\ Math.\ Oxford (2) 39 (1988) 185 - 199, Pacific J.\ Math.\ 134
(1988) 377 - 392}

\ref {Deu1}{Deundyak, V.M.}{Computation of the homotopy groups of the set of
invertible elements of certain $C^\a$-algebras of operators, and applications}
{Russ.\ Math.\ Surveys 35,3 (1980) 217 - 222}

\ref{Deu2}{Deundyak, V.M.}{Contractibility of some homotopically non trivial
groups of invertible operators with respect to groups containing them}{Algebra
and Discrete Math., pp.\ 46 - 54, Kalmytsk.\ Gos.\ Univ., Elisto 1985
(Russian)} 

\ref{DD}{Dixmier, J., Douady, A.}{Champs continus d'espaces hilbertiens et de
$C^\a$-alg\`ebres}{Bull.\ Soc.\ Math.\ France 91 (1963) 227 - 284}

\ref {Eck1}{Eckmann, B.}{Zur Homotopietheorie gefaserter R\"aume}{Comment.\
Math.\ Helv.\ 14 (1941/42) 141 - 192}

\ref {Eck2}{Eckmann, B.}{\"Uber die Homotopiegruppen der
Gruppenr\"aume}{Comment.\ Math.\ Helv.\ 14 (1941/42) 234 - 256}

\ref {Eck3}{Eckmann, B.} {Espaces fibr\'es et homotopie}{In: Colloque de
Topologie (espaces fibr\'es), Bruxelles 1950, 83 - 89, Masson et Cie, Paris,
1951}

\ref {Ekm}{Ekman, K.E.}{Unitaries and partial isometries in a real
$W^\a$-algebra}{Proc.\ Amer.\ Math.\ Soc.\ 54 (1976) 138 - 140}

\ref {Ell}{Elliott, E.G.}{On the classification of $C^\a$-algebras of real rank
zero}{Preprint}

\ref {ET}{Elworthy, K.D., Tromba, A.J.}{Differential structures and Fredholm
maps on Banach manifolds}{Global analysis, Berkeley 1968, Proc.\ Symp.\
Pure Math.\ 15, pp.\ 45 - 94, Amer.\ Math.\ Soc., Providence, R.I., 1970}

\ref {EFT}{Enomoto, M., Fuji, M., Takehana, H.}{On a conjecture of
Breuer}{Math.\ Japonicae 21 (1976) 387 - 389}

\ref {For}{Forster, O.}{Funktionentheoretische Hilfsmittel in der Theorie der
kommutativen Banachalgebren}{Jber.\ DMV 76 (1974) 1 - 17}

\ref {Fuj1}{Fujii K.}{Note on a paper of J.L.\ Taylor}{Mem.\ Fac.\ Sci.\ Kyushu Univ.\ Ser.\ A 32
(1978) 123 - 136}

\ref {Fuj2}{Fujii, K.}{A representation of complex $K$-groups by means of a
Banach algebra}{Mem.\ Sci.\ Kyushu Univ.\ Ser.\ A 32 (1978) 255 - 265}

\ref {Fuj3}{Fujii, K.}{A representation of real $K$-groups by means of a Banach
algebra}{Mem.\ Fac.\ Sci.\ Kyushu Univ.\ Ser.\ A 34 (1980) 379 - 394}

\ref {FN} {Furukawa, Y., Nomura, Y.}{Some homotopy groups of complex Stiefel
manifolds I - IV (III, IV without Furukawa)}{Sci.\ Rep.\ College Ed.\ Osaka Univ.\
25 (1976) 1 - 17, 27 (1978) 33 - 48, Bull.\ Aichi Ed.\ Natur.\ Sci.\ 31 (1982)
47 - 63, 32 (1983) 43 - 61}

\ref {Frt} {Furutani, K.}{A note on the Arens--Royden theorem for real Banach
algebras}{TRU Math.\ 11 (1975) 5 - 8}

\ref {Geb}{Geba, K.}{On the homotopy groups of $GL_c(E)$}{Bull.\ Acad.\ Polen.\
Sci.\ Math.\ 16 (1968) 699 - 702 (see also: Fund.
Math. 64 (1969) 341 - 373)}

\ref {Gew}{Gewalter, L.}{Der Periodizit\"atssatz f\"ur
$AW^\a$-Algebren}{Dissertation, Univ.\ Marburg, 1982}

\ref{Gio}  {Giordano, T.} {A classification of approximately finite real
$C^\a$-algebras}{J.\ reine angew.\ Math.\ 385 (1988) 161 - 194 }

\ref{GK}{Gohberg, I.C., Krupnik, N.Ya.}{Einf\"uhrung in die Theorie der
eindimensionalen singul\"aren Integraloperatoren}{Birkh\"auser, Basel, 1979}

\ref {Gra}{Grachev, V.A.}{Connectivity of the group of automorphisms of a
nuclear Fr\'echet space with a basis}{Math.\ Notes 35 (1984) 272 - 278}

\ref {Grm}{Gramsch, B.}{Relative Inversion in der St\"orungstheorie von
Operatoren und $\psi$-Algebren}{Math.\ Ann.\ (1984) 27 - 71}

\ref {Han1}{Handelman, D.E.} {$K_0$ of von Neumann algebras and $AF \,
C^\a$-algebras}{Quart.\ J.\ Math.\ Oxford (2) 29 (1978) 429 - 441}

\ref {Han2}{Handelman, D.E.}{Stable range in $AW^\a$-algebras}{Proc.\ Amer.\
Math.\ Soc.\ 76 (1979) 241 - 249}

\ref {dlH1}{de la Harpe, P.}{Classical Banach--Lie Algebras and Banach--Lie
Groups}{Lect.\ Notes Math.\ 285, Springer-Verlag, Berlin, 1972}

\ref {dlH2}{de la Harpe, P.}{The Clifford algebra and the spinor group
of a Hilbert space}{Composito Math.\ 25 (1972) 245 - 261}

\ref {HK}{Harris, L.A., Kaup, W.}{Linear algebraic groups in infinite
dimensions}{Illinois J.\ Math.\ 21 (1977) 666 - 674}

\ref {Hcz}{Hurewicz, W.}{Beitr\"age zur Topologie der Deformationen I -
IV}{Nederl.\
Akad.\ We\-tensch.\ Proc.\ Ser.\ A 38 (1935) 112 - 119, 521 - 528, 39 (1936) 117
- 126, 215 - 224}

\ref {Htz}{Hurwitz, A.}{\"Uber die Erzeugung der Invarianten durch
Integration}{Nachr.\ G\"ott.\ Ges.\ Wiss., Math.\ - Phys.\ Kl.\ 1897, 71 - 90,
also in: Mathem. Werke, Bd 2, LXXXI}

\ref {Ing}{Ingelstam, L.}{Real Banach algebras}{Ark.\ Math.\ 5 (1964) 239 - 279}

\ref {Jam}{James, I.M.}{Note on factor spaces}{J.\ London Math.\ Soc.\ 28
(1953) 278 - 285}

\ref {KK}{Kamei, E., Kato, Y.}{Homotopical properties of partial isometries of
the Calkin algebra}{Math.\ Japonicae 22 (1977) 83 - 87}

\ref {Kan}{Kandelaki, T.K.}{$K$-theory of ${\Bbb Z}_2$-graded Banach
categories I}{In: K-Theory and Homological Algebra, Tbilisi 1987-88, Lect.\
Notes Math.\ 1437, pp.\ 180 - 221, Springer-Verlag, Berlin, 1990}

\ref {Kar1}{Karoubi, M.}{Alg\`ebres de Clifford et $K$-th\'eorie}{Ann.\ Sci.\
Ecole Norm.\ Sup.\ (4) 1 (1968) 161 - 270}

\ref {Kar2}{Karoubi, M.}{Espaces classifiants en $K$-th\'eorie}{Trans.\ Amer.\
Math.\ Soc.\ 147 (1970) 74 - 115}

\ref {Kar3}{Karoubi, M.}{$K$-theory: An introduction}{Springer-Verlag, New
York--Berlin--Heidelberg 1978}

\ref {Kar4}{Karoubi, M.}{K-th\'eorie alg\'ebrique de certaines alg\`ebres
d'op\'erateurs}{In: Alg\`ebres d'Op\'erateurs, Les Plans-sur-Bex, 1978, Lect.\
Notes Math.\ 725, pp.\ 254 - 290, Springer-Verlag, Berlin, 1979}

\ref {Ksp1}{Kasparov, G.G.}{Lorentz groups, $K$-theory of unitary
representations and crossed products}{Dokl.\ Akad.\ Nauk SSSR 275 (1984) 541 -
545}

\ref{Ksp2} {Kasparov, G.G.}{Equivariant $KK$-theory and the Novikov
conjecture}{Invent.\ Math.\ 41 (1988) 147 - 201}

\ref{KSk}{Kasparov, G.G., Skandalis, G.}{Groups acting on buildings, operator
K-theory, and Novikov conjecture}{K-Theory 4 (1990/91) 303 - 337}

\ref{Khi1}{Khimshiashvili, G.N.}{Polysingular operators and the topology of
invertible singular operators}{Zeit.\ Anal.\ Anwend.\ 5,2 (1986) 139 - 145}

\ref {Khi2}{Khimshiashvili, G.N.}{On the topology of invertible linear singular
integral operators}{In: Global Analysis -- Studies and Appl.\ II, Lect.\ Notes
Math.\ 1214 (1987) 211 - 236}

\ref {Ku\v c}{Ku\v cment, P.A.}{A remark on the homotopy type of the
group of J-unitary operators}{Mat.\ Issled.\ 9, 4 (34) (1974) 170 -
171 (Russ.)}

\ref {KP}{Ku\v cment,
P.A., Pankov, A.A.}{Classifying spaces for equivariant $K$-theory}{Math.\ USSR
Sbornik 24 (1974) 31 - 48 (Letter to the editor, 27 (1975) 564)}

\ref {Kui}{Kuiper, N.H.}{The homotopy type of the unitary group of a Hilbert
space}{Topology 3 (1965) 19 - 30}

\ref {LZ}{Larotonda, A.R., Zalduendo, I.}{Spectral sets as Banach
manifolds}{Pac.\ J.\ Math.\ 120 (1985) 401 - 416}

\ref {Les}{Lesch, M.}{K-theory of Toeplitz $C^{\a}$-algebras on Lie
spheres}{Int.\ Equat.\ Oper.\ Theory 14 (1991) 120 - 145}

\ref {Luf1}{Luft, E.}{Maximal $R$-sets, Grassmann spaces, and Stiefel spaces of
a Hilbert space}{Trans.\ Amer.\ Math.\ Soc.\ 126 (1967) 73 - 107}

\ref{Luf2}{Luft, E.}{On the structure of maximal $R$-sets of a Hilbert
space}{Math.\ Ann.\ 175 (1968) 220 - 238}

\ref {Lun}{Lundell, A.T.}{Concise tables of James numbers and some
homotopy of classical Lie groups and associated homogeneous
spaces}{Lect. Notes Math. 1509, pp. 250 - 272, Springer-Verlag, Berlin,
1992}

\ref{Mat} {Matsumoto, K.}{Non-commutative three dimensional spheres}
{Japan. J. Math. 17 (1991) 333 - 356}

\ref{MT} {Matsumoto, K., Tomiyama,J.}{Non-commutative lens spaces}
{J. Math. Soc. Japan 44 (1992) 13 - 41}

\ref {Min}{Mingo, J.A.}{$K$-theory and multipliers of stable
$C^\a$-algebras}{Trans.\ Amer.\ Math.\ Soc.\ 299 (1987) 397 - 411}

\ref {Mit}{Mitjagin, B.S.}{The homotopy structure of the linear group of a
Banach space}{Russ.\ Math.\ Surveys 25 (1970) 59 - 103}

\ref {Nis}{Nistor, V.}{On the homotopy groups of the automorphism group of
$AF\, C^{\ast}$-algebras}{J.\ Operator Theory 19 (1988) 319 - 340}

\ref {Nom}{Nomura, Y.}{Some homotopy groups of real Stiefel manifolds in the
metastable range I - VI}{Sci.\ Rep.\ College Gen.\ Ed.\ Osaka Univ.\ 27 (1978)
1 - 31, 55 - 97, 28 (1979) 1 - 26, 35 - 60, 29 (1980) 159 - 183, 30 (1981) 11 -
57}

\ref {\^Ogu}{\^Oguchi, K.}{Homotopy groups of $Sp(n)/Sp(n-2)$}{J.\ Fac.\ Sci.\
Univ.\ Tokyo 16 (1969) 179 - 201}

\ref {Pal}{Palais, R.S.}{On the homotopy type of certain groups of
operators}{Topology 3 (1965) 271 - 279}

\ref {Pau}{Paulsen, V.I.}{The group of invertible elements in a Banach
algebra}{Coll.\ Math.\ 47 (1982) 97 - 100}

\ref{Php}{Phillips, J.}{K-theory
relative to a semifinite factor}{Indiana Univ.\ Math.\ J.\ 39 (1990) 339 - 354}

\ref {Phs1}{Phillips, N.C.}{K-theory of Fr\'echet algebras}{Int.\ J.\ Math.\ 2
(1991) 77 - 129}

\ref {Phs2}{Phillips, N.C.}{Five problems on operator algebras}
{Contemp.\ Math.\ 120, pp.\ 133 - 138, Amer.\ Math.\ Soc., Providence, R.I.,
1990}

\ref {Phi}{Phillips, R.S.}{On symplectic mappings of contraction
operators}{Studia Math.\ 31 (1968) 15 - 27}

\ref{Pim1}{Pimsner, M.V.} {Ranges of traces on $K_0$ of reduced crossed
products by free groups}{In: Operator Algebras and their Applications to Topology and
Ergodic Theory, Busenti, 1983, Lect. Notes Math. 1132, pp.\ 374 - 408,
Springer-Verlag, Berlin, 1985}

\ref{Pim2}{Pimsner, M.V.} { KK-groups of crossed products by groups
acting on trees}{Invent.\ Math.\ 86 (1986) 603 - 634}

\ref {PV1}{Pimsner, M.V., Voiculescu, D.}{Exact
sequences for $K$-groups and Ext-groups of certain cross-product
$C^\a$-algebras}{J.\ Operator Theory 4 (1980) 93 - 118}

\ref {PV2}{Pimsner, M.V., Voiculescu, D.}{$K$-groups of reduced crossed
products by free groups} {J.\ Operator Theory 8 (1982) 131 - 156}

\ref{Ply}{Plymen, R.J.}{Some recent results on infinite-dimensional spin
groups}{In: Adv.\ Math.\ Suppl.\ Studies 6, 159 - 171, Academic Press, New
York, 1979}

\ref{PR1}{Porta, H., Recht, L.}{Spaces of projections in a Banach
algebra}{Acta Cient.\ Venezolana 38.4 (1987)
408 - 426}

\ref{PR2}{Porta, H., Recht, L.}{Continuous selections of complementary
subspaces}{In: Contemp.\ Math.\ 52, 121 - 125, Amer.\ Math.\ Soc., Providence,
R.I., 1986}

\ref {PrS}{Pressley, A., Segal, G.B.}{Loop groups}{Oxford Math.\ Mono., Oxford
Sci.\ Publ., Claredon Press, Oxford 1986}

\ref {PW}{Putnam, C.R., Wintner, A.}{The connectedness of the orthogonal group
in Hilbert space}{Proc.\ Natl.\ Acad.\ Sci.\ USA 37 (1951) 110 - 112}

\ref {Rae1}{Raeburn, I.}{The relationship between a commutative Banach algebra
and its maximal ideal space}{J.\ Funct.\ Anal.\ 25 (1977) 366 - 390}

\ref {Rae2}{Raeburn, I.}{$K$-theory and $K$-homology relative to a
II$_\infty$-factor}{Proc.\ Amer.\ Math.\ Soc.\ 71 (1978) 294 - 298}

\ref {Ric}{Rickart, C.E.}{General theory of Banach algebras}{Van Nostrand,
Princeton 1960}
 
\ref {Rie1}{Rieffel,
M.A.}{Dimension and stable rank in the $K$-theory of $C^\a$-algebras}{Proc.\
London Math.\ Soc.\ (3) 46 (1983) 301 - 333}

\ref{Rie2}  { Rieffel, M.A.} { The homotopy groups of the unitary groups of
non-com\-mu\-tative tori}{ J.\ Operator Theory 17 (1987) 237 -
254}

\ref{Ros1}{Rosenberg, J.M.}{Homological invariants of extensions of
$C^\a$-algebras}{In:
Operator Algebras and Applications, Part I, Kingston 1980, Proc. Sympos.
Pure Math.\ 38,1, pp.\ 35 - 76,
Amer. Math. Soc., Providence, R.I., 1982}

\ref{Ros2}{Rosenberg, J.M.}{Group $C$*-algebras and topological invariants}{In:
Operator algebras and group representation II, Neptun 1980, pp.\ 95 - 115, 
Pitman, London, 1984}

\ref{Ros3}{Rosenberg, J.M.}{$C^\a$-algebras, positive scalar curvature, and
the Novikov conjecture III}{Topology 25 (1986) 319 - 336}

\ref{Rou}{Rouhani, A.}{Quasi-rotation $C^\a$-algebras}{Pac.\ J.\ Math.\ 148
(1991) 131 - 151}

\ref{Srb}{Schreiber, W.}{$K$-Theorie bez\"uglich Reeller von Neumannscher
Algebren}{Dissertation, Univ.\ Marburg, 1977}

\ref{Srd1} {Schr\"oder, H.}{On the homotopy type of the regular group of
a $W^\a$-algebra}{Math.\ Ann.\ 267 (1984) 271 - 277}

\ref{Srd2} {Schr\"oder, H.} {On the homotopy type of the regular group of
a real $W^\a$-algebra}{Int.\ Equat.\ Oper.\ Theory 9 (1986) 694 - 705}

\ref{Srd3} {Schr\"oder, H.} {On the topology of classical groups and
homogeneous spaces associated with a $W^\a$-algebra factor}{Int.\ Equat.\ Oper.\
Theory 10 (1987) 812 - 818}

\ref{Srd4} {Schr\"oder, H.} {K-Theory for Real $C^\a$-algebras and
Applications}{Pitman Research Notes Series 290, Longman, Harlow, 1993}

\reff{Srd5} {Schr\"oder, H.}{(unpublished)}

\ref{Sem}{Semenov, P.V.}{Contractibility of the linear group of spaces of
continuous functions on ordered compact sets}{In: Studies in the theory of
functions of several real variables, pp.\ 114 - 126, Yaroslav Gos.\ Univ.\
Yaroslavl', 1984 (Russ.)}

\ref{Sty}{Stacey, J.P.}{Stability of involutory $\ast$-antiautomorphisms
in
UHF-algebras}{J.\ Operator Theory 24 (1990) 57 - 74}

\ref{Ste}{Steenrod, N.E.}{The
topology of fibre bundles}{Princeton Univ.\ Press, Princeton 1951}

\ref{Str}{Stern, J.}{Le groupe des isometries d'un espace de Banach}{Studia
Math.\ 64 (1979) 139 - 149}

\ref{SS}{Sukochev, F.A., Sheremetyev, V.E.}{The linear groups of injective 
factors and of matroid $C^\a$-algebras are contractible to a point}{Math.\ 
Scand.\ 77 (1995) 119-128}

\ref{Swa}{Swanson, R.C.}{Linear symplectic structures on Banach spaces}{Rocky
Mountain J.\ Math.\ 10 (1980) 305 - 317}

\ref{Tay1}{Taylor,J.L.}{Banach algebras and topology}{In: Adv.\ Math.\ Studies, Algebras in
Analysis, 118 - 186, Academic Press, New York, 1975}

\ref{Tay2}{Taylor, J.L.}{Topological invariants of the maximal ideal space of a
Banach algebra}{Adv.\ Math.\ 19 (1976) 149 - 206}

\ref{Tay3}{Taylor, J.L.}{Twisted products in Banach algebras and third \v Cech
cohomology}{In: Lect.\ Notes Math.\ 575, pp.\ 157 - 174, Springer-Verlag, 
Berlin, 1977}

\ref{Tho}{Thomsen, K.}{Non-stable K-theory for operator algebras}
{K-Theory 4 (1991) 245 - 267}

\ref{Tod}{Toda, H.}{Quelques tables des groupes d'homotopie des groupes de
Lie}{C.R.\ Acad.\ Sci.\ Paris 241 (1955) 922 - 923}

\ref{Tro}{Troitsky, E.V.}{Kuiper's theorem for Hilbert modules: the general 
case}{Preprint MPI 96-16, Max-Planck-Institut Bonn, 1996}

\ref{Was}{Wassermann, A.}{A proof of the Connes-Kasparov conjecture for
connected linear Lie groups}{C.R.\ Acad.\ Sci.\ Paris S\'er.\ I Math.\ 304
(1987) 559 - 562}

\ref{Wey1}{Weyl, H.}{Theorie der Darstellung kontinuierlicher halbeinfacher
Gruppen durch lineare Transformationen I - III}{Math.\ Z.\ 23 (1925) 271 - 309,
24 (1926) 328 - 395, also in: Gesammelte Abhandlungen II, pp.\ 543 - 647, see
also pp. 453 - 467 for the earlier announcements}

\ref{Wey2}{Weyl, H.}{The Classical Groups}{Princeton Univ.\ Press, Princeton
1946}

\ref{Whi}{ Whitehead, G.W.}{Homotopy properties of the real orthogonal
groups}{Ann.\ of Math.\ 43 (1942) 132 - 146}

\ref{Wil}{Willgerodt, W.}{\"Uber den Homotopietyp der Automorphismengruppe
einer $W^\a$-Algebra}{Dissertation, Univ.\ Marburg, 1977}

\ref{Win1}{Wintner, A.}{Zur Theorie der beschr\"ankten Bilinearformen}{Math.\
Z.\ 30 (1929) 228 - 284}

\ref{Win2}{Wintner, A.}{\"Uber die automorphen Transformationen beschr\"ankter
nicht-singul\"arer hermitescher Formen}{Math.\ Z.\ 38 (1934) 695 - 700}

\ref{Win3}{Wintner, A.}{On bounded skew-symmetric forms}{Proc.\ Edinburgh
Math.\ Soc.\ II, 5 (1937) 90 - 92}

\ref{Woj}{Wojciechowski, K.}{A note on the space of pseudodifferential
projections with the same principal symbol}{J.\ Operator Theory 15 (1986) 207 -
216}

\ref{Woo} {Wood, R.} {Banach algebras and Bott periodicity}{ Topology 4
(1966) 371 - 389}

\ref{Yue}{Yuen, Y.}{Group of invertible elements of Banach algebras}{Bull.\
Amer.\ Math.\ Soc.\ 79 (1973) 82 - 84}

\ref{ZKKP}{Zaidenberg, M.G., Krein, S.G., Ku\v cment, P.A., Pankov, A.A.}
{Banach bundles and linear operators}{Russ.\ Math.\ Surveys 30 (1975) 115 - 175}

\ref{Zha1}{Zhang, S.}{Certain $C^\ast$-algebras with real rank zero and their
corona and multiplier algebras I,II}{Pac. J. Math. 155 (1992) 169 - 197, 
$K$-Theory 6 (1992) 1 - 27}

\ref{Zha2}{Zhang, S.}{On the homotopy type of the unitary group and the
Grassmann space of purely infinite simple $C^\a$-algebras}{$K$-Theory 
(to appear)}

\ref{Zha3}{Zhang, S.}{Matricial structure and homotopy type of simple
$C^\ast$-algebras with real rank zero}{J.\ Operator Theory 26 (1991) 283 - 312}

\ref{Zha4}{Zhang, S.}{K-theory and homotopy of certain groups and infinite
Grassmann spaces associated with $C^\a$-algebras}{Intern.\ J.\ Math.\ 5 (1994) 
425 - 445}

\bigskip\vskip 1cm

\leftline{Author's Adress:}
\leftline{Dr.\ Herbert Schr"der}
\leftline{c/o Fachbereich Mathematik der Universit\"at Dortmund}
\leftline{Vogelpothsweg 87, D - 44221 Dortmund}
\leftline{e-mail: schroed\@mathematik.uni-dortmund.de}

\end